\documentclass[11pt,fleqn,twoside]{article}
\usepackage{amsfonts,amssymb,latexsym}
\makeatletter
\newcommand{\prava}[1]{\small\it
\begin{flushleft}
Copyright \copyright \ 1999 by  #1
\end{flushleft}}

\newcommand{\name}[1]{\begin{flushleft}
                       \LARGE \bf #1
                       \end{flushleft}\vspace{-3mm}}

\newcommand{\Author}[1]{\begin{flushleft}
                       \it #1 \end{flushleft}}

\newcommand{\Adress}[1]{\begin{flushleft}
                       \it #1 \end{flushleft}}

\newcommand{\Date}[1]{\begin{flushleft}
                      \small  \it #1 \end{flushleft}}

\newcommand{\ehkol}{Author \ name}
\newcommand{\ohkol}{Article \ name}
\renewcommand{\@evenhead}{
\hspace*{-3pt}\raisebox{-15pt}[\headheight][0pt]{\vbox{\hbox to \textwidth 
{\thepage \hfil \ehkol}\vskip4pt \hrule}}}
\renewcommand{\@oddhead}{
\hspace*{-3pt}\raisebox{-15pt}[\headheight][0pt]{\vbox{\hbox to \textwidth 
{\ohkol \hfil \thepage}\vskip4pt\hrule}}}
\renewcommand{\@evenfoot}{}
\renewcommand{\@oddfoot}{}

\newcommand{\be}{\begin{equation}}
\newcommand{\ee}{\end{equation}}
\newcommand{\ba}{\hspace*{-5pt}\begin{array}}
\newcommand{\ea}{\end{array}}

\newcommand{\ds}{\displaystyle}
\makeatother

\newcommand{\pbf}[1]{\mbox{\mathversion{bold}$#1$}}

\begin{document}

\thispagestyle{empty}
\setcounter{page}{448}
\renewcommand{\ehkol}{B.A. Kupershmidt}
\renewcommand{\ohkol}{What a Classical $r$-Matrix Really Is}

\begin{flushleft}
\footnotesize \sf
Journal of Nonlinear Mathematical Physics \qquad 1999, V.6, N~4,
\pageref{kupershmidt_7-fp}--\pageref{kupershmidt_7-lp}.
\hfill {\sc Article}
\end{flushleft}

\vspace{-5mm}

\renewcommand{\footnoterule}{}
{\renewcommand{\thefootnote}{}
 \footnote{\prava{B.A. Kupershmidt}}}

\name{What a Classical {\mathversion{bold}$r$}-Matrix Really Is}\label{kupershmidt_7-fp}

\Author{Boris A. KUPERSHMIDT}

\Adress{Department of Mathematics, University of Tennessee
Space Institute,\\ 
Tullahoma, TN 37388,  USA\\
E-mail: bkupersh@utsi.edu}

\Date{Received September 08, 1999; Accepted October 14, 1999}

\begin{flushright}
\begin{minipage}{7cm}
\small  \bfseries \itshape
To my friend and colleague
K.C.~Reddy on occasion of
his retirement.
\end{minipage}
\end{flushright}

\begin{abstract}
\noindent
The notion of classical $r$-matrix is re-examined, and a def\/inition suitable to 
dif\/fe\-ren\-tial (-dif\/ference) Lie algebras, -- where the standard def\/initions are shown to 
be def\/icient,~-- is proposed, the notion of an ${\mathcal O}$-operator.  This notion has all 
the natural properties one would expect form it, but lacks those which are artifacts 
of f\/inite-dimensional isomorpisms such as not true in dif\/ferential generality 
relation $\mbox{End}\, (V) \simeq V^* \otimes V$ for a vector space $V$.  
Examples considered include a quadratic Poisson bracket on the dual space to a 
Lie algebra;     generalized symplectic-quadratic models of 
such brackets 
(aka Clebsch representations); and Drinfel'd's 2-cocycle interpretation 
of nondegenate classical $r$-matrices.
\end{abstract}

\centerline{\bf Contents}
\begin{enumerate}
\item[\S~1.] Form Artin relation to Quantum Yang--Baxter equation to Classical 
Yang--Baxter equation
\item[\S~2.] Classical $r$-matrices and 2-cocycles
\item[\S~3.] Dif\/ferential Lie algebras say: 2-cocycles, -- Si, $r$-matrices, -- No, 
${\mathcal O}$-operators are welcome
\item[\S~4.] ${\mathcal O}$-natural property of the ${\mathcal O}$-operators
\item[\S~5.] Linear Poisson  brackets on dual spaces to Lie algebras
\item[\S~6.] Quadratic Poisson brackets on dual spaces to Lie algebras
\item[\S~7.] Symplectic models for linear Poisson brackets on dual spaces to Lie algebras
\item[\S~8.] Clebsch representations for quadratic Poisson brackets on dual spaces to Lie algebras
\item[\S~9.] Properties of the quadratic Poisson brackets on $V \oplus V^*$
\item[]\hspace*{-8mm}Appendix A1. Crossed Lie algebras
\item[]\hspace*{-8mm}Appendix A2. Symplectic $r$-matrices and symplectic doubles
\end{enumerate}

\renewcommand{\theequation}{\arabic{section}.\arabic{equation}}
\setcounter{equation}{0}

\section{From Artin relation to Quantum Yang--Baxter equation\\
\hspace*{12mm}to Classical Yang--Baxter equation}

\rightline{Hier ist kein Warum.}

\medskip

This paper is written with a non-expert in mind, and the text is purposedly self-contained 
apart from a few references to basic properties of the algebraic calculus of 
variations and Hamiltonian formalism.  In this section we derive the Classical 
Yang--Baxter equation (CYBE) as the quasiclassical limit of the Quantum Yang-Baxter equation 
(QYBE); the latter will be seen in a moment as a special form of the Artin relation for 
the generators of the braid group.

We thus start from a purely f\/inite-dimensional view-point; f\/inite-dimensional Lie 
algebras will come in later on through an interpretation of the CYBE, and 
dif\/ferential Lie algebras will appear later still.  Have I mentioned that most, 
if not all, of the results in this Section are gleamed from the conf\/idential list 
of examination questions given annually to all low-level NSA employees?

Let's f\/ix a vector space $V$ and let
\be
S:\ \ V \otimes V \rightarrow V \otimes V 
\ee
be an operator.  $S$ induces the operators $S^{12} = S \otimes  {\mathbf 1}$ and 
$S^{23} = {\mathbf 1} \otimes S$ acting on $V \otimes V \otimes V$ in an obvious way. 
 The operator $S^{13}$ acts on $V \otimes V \otimes V$ in an equally natural way:
\be
S^{13} (e_{i} \otimes e_{j} \otimes e_{k}) = \sum_{ab} S^{ab}_{ik} 
e_{a} \otimes e_{j} \otimes e_{b}, 
\ee
where
\be
S (e_{i} \otimes e_{j}) = \sum_{cd} S^{cd}_{ij} e_{c} \otimes e_{d}, 
\ee
and $(e_{i})$ is a basis in $V$.

Denote by $P: V \otimes V \rightarrow V \otimes V$ the permutation operator,
\be
P (v_{1} \otimes v_{2}) = v_{2} \otimes v_{1}, \qquad \forall \; v_{1}, v_2 
\in V, 
\ee
and let ${\mathcal M}: V \otimes V \otimes V \rightarrow V \otimes V \otimes V$
be the operator of mirror symmetry:
\be
{\mathcal M} (v_1 \otimes v_2 \otimes v_3) = v_3 \otimes v_2 \otimes v_1. 
\ee

It's immediate to see by inspection that
\be
{\mathcal M} = P^{23} P^{12} P^{23} = P^{12} P^{23} P^{12}. 
\ee
The Artin equation for an arbitrary operator
$S: V \otimes V \rightarrow V \otimes V$ is then
\be
 S^{23} S^{12} S^{23} = S^{12} S^{23} S^{12}. 
\ee

That's all there is to it.  We now proceed to massage this equation in various directions.

Set
\be
S = P R \Leftrightarrow R = P S. 
\ee
The Artin equation (1.7) will become then
\be
P^{23} R^{23} P^{12} R^{12} P^{23} R^{23} = P^{12} R^{12} P^{23} R^{23} P^{12} R^{12}. 
\ee

Let $A$, ${\mathcal U}$, ${\mathcal V}$, ${\mathcal W}$ be arbitrary operators $V \otimes V
 \rightarrow V \otimes V$.  The following identities are easy to check and are left 
for the reader to verify:
\renewcommand{\theequation}{\arabic{section}.\arabic{equation}{\rm a}}
\setcounter{equation}{9}
\be
A^{12} P^{23} = P^{23} A^{13}, 
\ee
\renewcommand{\theequation}{\arabic{section}.\arabic{equation}{\rm b}}
\setcounter{equation}{9}
\be
A^{13} P^{23} = P^{23} A^{12},
\ee
\renewcommand{\theequation}{\arabic{section}.\arabic{equation}{\rm c}}
\setcounter{equation}{9}
\be
A^{23} P^{12} = P^{12} A^{13},
\ee
\renewcommand{\theequation}{\arabic{section}.\arabic{equation}{\rm d}}
\setcounter{equation}{9}
\be
A^{13} P^{12} = P^{12} A^{23}; 
\ee
\renewcommand{\theequation}{\arabic{section}.\arabic{equation}}
\setcounter{equation}{10}
\be
{\mathcal U}^{23} {\mathcal V}^{12} {\mathcal W}^{23} = {\mathcal W}^{12} 
{\mathcal V}^{23} {\mathcal U}^{12}
\ee
whenever at least two out of three operators ${\mathcal U}$, ${\mathcal V}$, 
${\mathcal W}$  are equal to $P$.

Notice that each one of the formulae (1.10) can be taken as an 
{\it invariant definition} of the operator $A^{13}$, thus avoiding the 
coordinate def\/inition (1.2).

Now let us transform separately each side of the equation (1.9).  For the LHS we get
\renewcommand{\theequation}{\arabic{section}.\arabic{equation}{\rm L}}
\setcounter{equation}{11}
\be
\ba{l}
P^{23} R^{23} P^{12} R^{12} P^{23} R^{23} \  {\mathop{=}\limits^{\mbox{\scriptsize [by (1.10c,a)]}}}
\  P^{23} P^{12} R^{13} P^{23} R^{13} R^{23} 
\vspace{2mm}\\
\ds \qquad \ \  {\mathop{=}\limits^{\mbox{\scriptsize [by  (1.10b)]}}} \
 P^{23} P^{12} P^{23} R^{12} R^{13} R^{23} = {\mathcal M} R^{12} R^{13} R^{23}, 
\ea
\ee
while for the RHS of the equation (1.9) we obtain
\renewcommand{\theequation}{\arabic{section}.\arabic{equation}{\rm R}}
\setcounter{equation}{11}
\be
\ba{l}
P^{12} R^{12} P^{23} R^{23} P^{12} R^{12} \   {\mathop{=}\limits^{\mbox{\scriptsize [by  (1.10a,c)]}}}
 P^{12} P^{23} R^{13} P^{12} R^{13} R^{12} 
\vspace{2mm}\\
\  {\mathop{=}\limits^{\mbox{\scriptsize [by  (1.10d)] }}}
\  P^{12} P^{23} P^{12} R^{23} R^{13} R^{12} =  {\mathcal M} R^{23} R^{13} R^{13}. 
\ea
\ee
Equating the expressions (1.12L) and (1.12R) we f\/ind
\renewcommand{\theequation}{\arabic{section}.\arabic{equation}}
\setcounter{equation}{12}
\be
R^{12} R^{13} R^{23} = R^{23} R^{13} R^{12}. 
\ee

This is called QYBE.  Since Artin equation (1.7) is satisf\/ie by $S=P$, the QYBE equation 
(1.13) is satisf\/ied by $R= {\mathbf 1}$ .  Let's look for {\it perturbations} 
of this solution: set
\be
R = {\mathbf 1} + h r + h^2 \rho + O \left(h^3\right) 
\ee
with some formal parameter $h$.  Then
\renewcommand{\theequation}{\arabic{section}.\arabic{equation}{\rm L}}
\setcounter{equation}{14}
\be
\ba{l}
\ds R^{12} R^{13} R^{23} = \left({\mathbf 1} + h r^{12} + h^2 \rho^{12}\right) 
\left({\mathbf 1} + hr^{13} + h^2 \rho^{13}\right) 
\left({\mathbf 1} + hr^{23} + h^2 \rho^{23}\right) +  O \left(h^3\right) 
\vspace{2mm}\\
\ds \qquad = {\mathbf 1} + h \left(r^{12} + r^{13} + r^{23}\right) + h^2 
\left(\rho^{12} + \rho^{13} + \rho^{23}\right) 
\vspace{2mm}\\
\ds \qquad + h^2 \left(r^{12} r^{13} + r^{12} r^{23} + r^{13} r^{23}\right) +O \left(h^3\right), 
\ea\hspace{-23.95pt}
\ee
\renewcommand{\theequation}{\arabic{section}.\arabic{equation}{\rm R}}
\setcounter{equation}{14}
\be
\ba{l}
\ds R^{23} R^{13} R^{12} = \left({\mathbf 1} + h r^{23} + h^2 \rho^{23}\right) 
\left({\mathbf 1} + h r^{13} + h^2 \rho^{13}\right) 
\left({\mathbf 1} + h r^{12} + h^2 \rho^{12}\right) + O \left(h^3\right)
\vspace{2mm}\\
\ds \qquad = {\mathbf 1} + h \left(r^{23} + r^{13} + r^{12}\right) + h^2 \left(\rho^{23} 
+ \rho^{13} + \rho^{12}\right) 
\vspace{2mm}\\
\ds \qquad + h^2 \left(r^{23} r^{13} + r^{23} r^{12} + r^{13} r^{12}\right) + O \left(h^3\right). 
\ea\hspace{-25.16pt}
\ee
Comparing the expressions (1.15L) and (1.15R), we see that they dif\/fer in $h^2$-terms; 
these yield
\renewcommand{\theequation}{\arabic{section}.\arabic{equation}}
\setcounter{equation}{15}
\be
c (r): = \left[r^{12}, r^{13}\right] + \left[r^{12}, r^{23}\right] + 
\left[r^{13}, r^{23}\right] = 0.
\ee
This is called CYBE.  As a quaslassical approximation to a noncommutative QYBE (1.13) 
in an associative framework, CYBE should have Poisson-brackets-related proprties and/or 
interpretations.  This is known to be true, and we shall see more of such presently.

The f\/irst step in this direction is to realize that since only {\it commutators}
are involved in the CYBE (1.16), the operator $r$, -- called a classical $r$-matrix, --
can be considered not just as an element of the tensor square of the Lie 
algebra $\mbox{End}\,(V)$:
\be
\mbox{End}\, (V) \otimes \mbox{End}\, (V) \approx \mbox{End}\, (V \otimes V), 
\ee
but as an element of ${\mathcal G} \otimes {\mathcal G}$ for {\it arbitrary} Lie algebra 
${\mathcal G}$:
\be
r \in {\mathcal G} \otimes {\mathcal G}. 
\ee
The CYBE (1.16) is then understood as an identity in ${\mathcal G} \otimes {\mathcal G}
 \otimes {\mathcal G}$:  If
\be
r = \sum_i a_i \otimes b_i, \qquad  a_i, b_i \in {\mathcal G}, 
\ee
then (temporarily stepping outside ${\mathcal G}$ into the
Universal enveloping algebra $U ({\mathcal G}))$
\renewcommand{\theequation}{\arabic{section}.\arabic{equation}{\rm a}}
\setcounter{equation}{19}
\be
\left[r^{12}, r^{13}\right] = 
\sum_{ij} [ a_i, a_j] \otimes b_i \otimes b_j, 
\ee
\renewcommand{\theequation}{\arabic{section}.\arabic{equation}{\rm b}}
\setcounter{equation}{19}
\be
\left[r^{12}, r^{23}\right] = \sum_{ij} a_i \otimes [b_i, a_j] \otimes b_j, 
\ee
\renewcommand{\theequation}{\arabic{section}.\arabic{equation}{\rm c}}
\setcounter{equation}{19}
\be
\left[r^{13}, r^{23}\right] = \sum_{ij} a_i \otimes a_j \otimes [b_i, b_j], 
\ee
so that the CYBE (1.16) becomes
\renewcommand{\theequation}{\arabic{section}.\arabic{equation}}
\setcounter{equation}{20}
\be
0 = c \left(\sum_i a_i \otimes b_i\right) 
= \sum_{ij} ([a_i, a_j] \otimes b_i \otimes b_j + a_i \otimes [b_i, a_j] 
\otimes a_j \otimes b_j + a_i \otimes a_j \otimes [b_i, b_j]). 
\ee
(Finding out the proper Lie-algebraic object to which the CYBE (1.16) 
in ${\mathcal G}^{\otimes 3}$ is the quasiclassical approximation is far from easy; 
this is done in Drinfel'd's paper~[3].)

\medskip

\noindent
{\bf Remark 1.22.}  Had we considered the quasiclassical approximation 
to the Artin equation (1.7) itself, in the form
\setcounter{equation}{22}
\be
S = P + h {\bar{r}} + h^2 {\bar{\rho}} + O \left(h^3\right),
\ee
the $h^2$-terms would have collected into the equation
\be
{\bar{r}}^{23} {\bar{r}}^{12} P^{23} + {\bar{r}}^{23} P^{12} 
{\bar{r}}^{23} + P^{23} {\bar{r}}^{12} {\bar{r}}^{23} =
 {\bar{r}}^{12} {\bar{r}}^{23} P^{12} + {\bar{r}}^{12} P^{23} {\bar{r}}^{12} 
+ P^{12} {\bar{r}}^{23} {\bar{r}}^{12}.
\ee
This equation turns into CYBE (1.16) upon the substitution
\be
{\bar{r}} = P r. 
\ee

Let us discuss the skewsymmetry proprty of the classical $r$-matrix.  
If we impose on the operator $S$ the very natural ``unitarity''condition
\be
S^2 = {\mathbf 1}, 
\ee
the $r$-matrix $\ds r = {\partial (PS) \over \partial h}\bigg|_{h=0}$ 
inherits from the unitarity the skewsymmetry condition
\be
P r = - r P. 
\ee
The Lie-algebraic $r$-matrix (1.19) then belongs to $\Lambda^2 {\mathcal G}$ rather than 
${\mathcal G}^{\otimes2}$:
\be
r = \sum_i a_i \wedge b_i = \sum_i (a_i \otimes b_i - b_i \otimes a_i). 
\ee
Although non-skewsymmetric $r$-matrices play many important r\^oles in 
various branches of Mathematics and Physics (see the textbook~[1] of Chari and Pressley 
as the basic reference for what follows), in this paper {\it all} $r$-matrices 
will be considered skewsymmetric, due to the nature of the topics discussed.

For future reference, we shall record the skewsymmetric version of formulae (1.20)
for $c (r)$ with the skewsymmetric $r$-matrix (1.28):
\renewcommand{\theequation}{\arabic{section}.\arabic{equation}{\rm a}}
\setcounter{equation}{28}
\be
\!\!\left[r^{12}, r^{13}\right] = \sum_{ij} ([a_i, a_j] \otimes b_i \otimes b_j - [a_i, b_j] 
\otimes b_i \otimes a_j - [b_i, a_j] \otimes a_i \otimes b_j + [b_i, b_j] \otimes a_
\otimes a_j), 
\ee
\renewcommand{\theequation}{\arabic{section}.\arabic{equation}{\rm b}}
\setcounter{equation}{28}
\be 
\!\!\left[r^{12}, r^{23}\right] = \sum_{ij} ( a_i \otimes [b_i, a_j] \otimes b_j - a_i \otimes 
[b_i, b_j] \otimes a_j - b_i \otimes [a_i, a_j] \otimes b_j + b_i \otimes [ a_i, b_j] \otimes a_j), 
\ee
\renewcommand{\theequation}{\arabic{section}.\arabic{equation}{\rm c}}
\setcounter{equation}{28}
\be
\!\!\left [r^{13}, r^{23}\right] = \sum_{ij} (a_i \otimes a_j \otimes [b_i, b_j] - a_i \otimes b_j 
\otimes [b_i, a_j] - b_i \otimes a_j \otimes [a_i, b_j] + b_i \otimes b_j \otimes 
[a_i, a_j]).  
\ee

\renewcommand{\theequation}{\arabic{section}.\arabic{equation}}
\setcounter{equation}{0}

\section{Classical {\mathversion{bold}$r$}-matrices and 2-cocycles}

Any skewsymmetric element $r \in \Lambda^2 {\mathcal G}$  
satisfying the CYBE (1.16) is called a classical $r$-matrix.

To every $r \in {\mathcal G} \otimes {\mathcal G}$ we can associate an operator 
${\mathcal O} = {\mathcal O}_r: {\mathcal G}^* \rightarrow {\mathcal G}$ by the rule
\be
\langle u, {\mathcal O} (v) \rangle = \langle u \otimes v, r \rangle, 
\qquad \forall  \; u, v \in {\mathcal G}^*. 
\ee
Conversely, this equality attaches an element $r \in {{\mathcal G}}^{\otimes 2}$ to every operator 
${\mathcal O}: {\mathcal G}^* \rightarrow {\mathcal G}$.  
(Why do such banalities deserve being mentioned?  Because 
they are not true in general.  Please bear with me).  The skewsymmetry of $r$ is equivalent to 
skewsymmetry of ${\mathcal O}$:
\be
\langle u, {\mathcal O} (v)\rangle + \langle v, {\mathcal O} (u) \rangle = 0. 
\ee

Now, suppose temporarily that $r$ is {\it nondegenerate}, i.e.,
 ${\mathcal O}$ is invertible.  (${\mathcal G}$ is then even-dimensional). 
 Consider the skewsymmetric bilinear form $\omega = \omega_r$ on~${\mathcal G}$:
\be
\omega (x, y) = \langle {\mathcal O}^{-1} (x), y \rangle, \qquad
 \forall \; x, y \in {\mathcal G}. 
\ee

\noindent
{\bf Theorem 2.4  (Drinfel'd [2]).}  {\it A nondegenerate $r$ in $\wedge^2 {\mathcal G}$ satisfies 
CYBE iff $\omega$ is a 2-cocycle on ${\mathcal G}$.}

\medskip

We postpone the proof of this Theorem until later on in this Section, since we aim at a higher 
prize:  to reformulate this 2-cocycle characterization of classical $r$-matrices into a form 
suitable for a fruitful def\/inition.

Let's write down the condition for $\omega$ to be a 2-cocycle on ${\mathcal G}$:
\setcounter{equation}{3}
\be
0 = \omega (x, [y, z)) + \mbox{c.p.} = \langle {\mathcal O}^{-1} (x), [y, z]\rangle + \mbox{c.p.}, 
\qquad \forall \; x, y, z \in {\mathcal G}. 
\ee
Since ${\mathcal O}$ is invertible, we can f\/ind $u, v, w \in {\mathcal G}^*$ such that 
\be
x = {\mathcal O} (u), \qquad y = {\mathcal O} (v), \qquad 
 z = {\mathcal O} (w). 
\ee
The 2-cocycle condition (2.4) then becomes
\be
\langle  u, [{\mathcal O} (v), {\mathcal O}(w)] \rangle + \mbox{c.p.} = 0,
\qquad  \forall \; u, v, w \in {\mathcal G}^*. 
\ee

What have we achieved?  First, in equality (2.6) the map ${\mathcal O}$ 
is no longer required to be invertible.  
(2-cocycles are often degenerate.  This is true in Fluid Mechanics and Plasma Physics, see 
many examples in~[6]; same happens in f\/inite dimensions, e.g. for complex semisimple Lie 
algebras, see discussion in~[1], p. 62.)  Second, the equaity (2.6) is trilinear in ${\mathcal G}^*$, 
but we can transform this trilinear equation into an equivalent bilinear one, as follows.

Denote the coadjoint action of $x \in {\mathcal G}$ on $u \in {\mathcal G}^*$ by $x^{.} u:$
\be
\langle x^{.} u, y\rangle = - \langle u, [x, y]\rangle, \qquad \forall \; y \in {\mathcal G}. 
\ee
Now, using the skewsymmetry of ${\mathcal O}$, we get
\renewcommand{\theequation}{\arabic{section}.\arabic{equation}{\rm a}}
\setcounter{equation}{7}
\be
 1) \quad  \langle u, [ {\mathcal O}(v), {\mathcal O}(w)]\rangle = 
- \langle {\mathcal O} (v)^{.} u, {\mathcal O} (w) \rangle = 
\langle w, {\mathcal O} ({\mathcal O}(v)^{.} u) \rangle, 
\ee
\renewcommand{\theequation}{\arabic{section}.\arabic{equation}{\rm b}}
\setcounter{equation}{7}
\be 
2) \quad \langle v, [{\mathcal O} (w), {\mathcal O}(u)]\rangle = 
\langle {\mathcal O} (u)^{.} v, {\mathcal O} (w) \rangle = 
- \langle w, {\mathcal O} ({\mathcal O}(u)^{.} v) \rangle. 
\ee
Substituting (2.8) into (2.6) we obtain
\renewcommand{\theequation}{\arabic{section}.\arabic{equation}}
\setcounter{equation}{8}
\be 
\langle w, {\mathcal O} ({\mathcal O} (v)^{.} u - {\mathcal O} (u)^{.} v) + 
[{\mathcal O} (u), {\mathcal O}(v) ] \rangle  = 0. 
\ee
Since $w$ is arbitrary, the 2-cocycle condition (2.6) is equivalent to the equality
\be
{\mathcal O} ({\mathcal O} (u)^{.} v - {\mathcal O} (v)^{.} u) = 
[ {\mathcal O} (u), {\mathcal O} (v) ], \qquad  \forall \; u, v \in {\mathcal G}^*. 
\ee

This suggests the following generalization of the notion of the classical $r$-matrix.  Let 
${\mathcal G}$ be a Lie algebra, ${\mathcal U}$ a ${\mathcal G}$-module, 
and ${\mathcal O}: {\mathcal U} \rightarrow {\mathcal G}$ a linear map.  Let's  
make ${\mathcal U}$ into an algebra by def\/ining a skew multiplication [~,~] in 
${\mathcal U}$ by the rule
\be
[u, v] = {\mathcal O} (u). v - {\mathcal O} (v). u, \qquad \forall \; u, v \in {\mathcal U}. 
\ee
${\mathcal O}$ is called an ${\mathcal O}$-operator, or a classical $r$-matrix, 
if\/f ${\mathcal O}$ is a homomorphism of algebras:
\be
{\mathcal O} ([u, v]) = [ {\mathcal O} (u), {\mathcal O} (v)], \qquad
 \forall \; u, v \in {\mathcal U}; 
\ee
equivalently, 
\be
{\mathcal O} ({\mathcal O} (u). v - {\mathcal O} (v). u) = 
[{\mathcal O} (u), {\mathcal O} (v)], \qquad \forall \; u, v \in {\mathcal U}. 
\ee

\noindent
{\bf Examle 2.14.} Let ${\mathcal G} = sl_2$ with a basis $(h; e; f)$:
\setcounter{equation}{14}
\be
[h, e] = 2e, \qquad  [h, f] = - 2f, \qquad [e, f] = h. 
\ee
Let ${\mathcal U}$ be 2-dimensional, with a basis $(v_0; v_1)$ and the action of 
$s\ell_2$ of the  fundamental representation:
\be
\ba{l}
e. v_0 = 0, \qquad h . v_0 = v_0, \qquad  f. v_0 = v_1, 
\vspace{2mm}\\
e. v_1 = v_0, \qquad  h . v_1 = - v_1, \qquad  f . v_1 = 0. 
\ea
\ee
Each of the following 2 maps can be easily seen to be an ${\mathcal O}$-operator:
\renewcommand{\theequation}{\arabic{section}.\arabic{equation}{\rm a}}
\setcounter{equation}{16}
\be
{\mathcal O}  \left(\matrix{ v_0 \cr v_1 \cr} \right) = c_1 \left(\matrix{h \cr f \cr} \right) + 
c_2 \left(\matrix{f \cr 0 \cr} \right), \qquad c_1, c_2 \  \mbox{are  constants}, 
\ee
\renewcommand{\theequation}{\arabic{section}.\arabic{equation}{\rm b}}
\setcounter{equation}{16}
\be
{\mathcal O}  \left(\matrix{ v_0 \cr v_1 \cr} \right) = c_3 \left(\matrix{e \cr -h \cr} \right) + 
c_4 \left(\matrix{0 \cr e \cr} \right),
 \qquad c_3, c_4 \  \mbox{are  constants}, 
\ee

\noindent
{\bf Proposition 2.18.} {\it If ${\mathcal O}: {\mathcal U} \rightarrow {\mathcal G}$ 
is an ${\mathcal O}$-operator then  ${\mathcal U}$ is a Lie algebra.}

\medskip

\noindent
{\bf Proof.}  By formula (2.11),
\renewcommand{\theequation}{\arabic{section}.\arabic{equation}}
\setcounter{equation}{18}
\[
\ba{l}
\ds [[u, v], w ] + \mbox{c.p.} = \{{\mathcal O} ([ u, v]). w - 
{\mathcal O} (w) . [u, v]\} + \mbox{c.p.} 
\vspace{2mm}\\
\ds \qquad {\mathop{=}\limits^{\mbox{\scriptsize [by (2.12)]}}}  \
\{ [ {\mathcal O} (u), {\mathcal O}(v)]. w - {\mathcal O} (w). 
({\mathcal O} (u). v - {\mathcal O} (v). u)\} + \mbox{c.p.}
\vspace{2mm}\\ 
\ds \qquad {\mathop{=}\limits^{\mbox{\scriptsize [since ${\mathcal U}$ is  
a  ${\mathcal G}$-module]}}} \   \{({\mathcal O} (u). ( {\mathcal O} (v). w) - 
{\mathcal O} (v). ( {\mathcal O}(u). w)) + \mbox{c.p.} \} 
\vspace{2mm}\\
\ds \qquad - \{ {\mathcal O}  (w). ({\mathcal O} (u). v) + \mbox{c.p.} \} 
+ \{ {\mathcal O} (w). ({\mathcal O} (v). u) + \mbox{c.p.} \}
\vspace{2mm}\\
\ds \qquad  = \{ ({\mathcal O}(u). ( {\mathcal O} (v). w) - 
{\mathcal O} (v) . ({\mathcal O} (u). w)) + \mbox{c.p.} \} 
\vspace{2mm}\\
\ds \qquad - \{ {\mathcal O}(u). ({\mathcal O}
(v). w) + \mbox{c.p.} \} + \{ {\mathcal O} (v). ({\mathcal O} (u). w) + \mbox{c.p.} \} = 0.
\hspace{96.5pt} \mbox{\rule{3mm}{3mm}}
\ea
\]

Thus, ${\mathcal O}$ is aposteriori a homomorphism of Lie algebras.  (This explains formulae
 (2.17), since 
${\mathcal U}$ is 2-dimensional and ${\mathcal O}$ 
maps ${\mathcal U}$ into $b_+$  or $b_-.$) We shall not wander into the 
general ${\mathcal U}$-route here.  (For example, adding {\it linear} 
on ${\mathcal O}$ conditions making ${\mathcal G} + {\mathcal U}$ 
into a Lie algebra.)  From now  on ${\mathcal U}$ reverts to the old familiar ${\mathcal G}^*$.  

\medskip

\noindent
{\bf Proposition 2.19.} {\it  Let ${\mathcal O}: {\mathcal G}^* \rightarrow {\mathcal G}$ be an ${\mathcal O}$-operator, so 
that ${\mathcal G}^*$ is now a Lie algebra.  Then the skewsymmetric bilinear form $\Omega$ on ${\mathcal G}^*$:
\setcounter{equation}{19}
\be
 \Omega (u, v) = \langle u, {\mathcal O} (v) \rangle 
\ee
is a 2-cocycle on ${\mathcal G}^*$.}

\medskip

\noindent
{\bf Proof.}  We have:
\[
\ba{l}
\ds \Omega ([u, v], w) + \mbox{c.p.} = \langle [u, v], {\mathcal O} (w) \rangle + 
\mbox{c.p.} = - \langle w, {\mathcal O} ([u, v]) \rangle  + 
\mbox{c.p.}
\vspace{2mm}\\
\ds \qquad {\mathop{=}\limits^{\mbox{\scriptsize [by (2.12)]}}} \
 - \langle  w,[{\mathcal O} (u), {\mathcal O}(v) ] \rangle  + \mbox{c.p.}
 \ {\mathop{=}\limits^{\mbox{\scriptsize [by (2.6) ]}}} \ 0. \hspace{146pt} \mbox{\rule{3mm}{3mm}}
\ea
\]

Let us prove now Drinfel'd's Theorem~2.4.  We shall evaluate each of the 3 terms in the 
${\mathcal O}$-equation (2.13) and compare them to the expressions (1.29) for $c(r$).  

First, for $r = \sum\limits_i (a_i \otimes b_i - b_i \otimes a_i)$ (1.28), we get
\[
\ba{l}
\ds \langle u \otimes v, r\rangle  = \sum_i (\langle  u, a_i \rangle \langle v, b_i\rangle 
- \langle u, b_i \rangle \langle  v, a_i \rangle ) 
\vspace{3mm}\\
\ds \qquad = \langle u, \sum_i (\langle v, b_i \rangle  a_i - \langle v, a_i \rangle b_i) \rangle, 
\ea
\]
so that
\be
{\mathcal O} (v) = \sum_i (\langle v, b_i \rangle a_i - \langle v, a_i \rangle b_i) .
\ee
Now, 
\renewcommand{\theequation}{\arabic{section}.\arabic{equation}{\rm a}}
\setcounter{equation}{21}
\be
\ba{ll}
1) & \ds  {\mathcal O} ({\mathcal O} (u)^\cdot v) = {\mathcal O} \left(\sum_j \left(\langle u, b_j \rangle
 (a_j \, ^{\cdot} v) - \langle u, a_j \rangle (b_j \, ^{\cdot}  v)\right)\right) 
\vspace{3mm}\\
& \ds \quad = \sum_{ji} (\langle u, b_j \rangle \langle a_j \, ^{\cdot} v,  b_i\rangle a_i - 
\langle u, a_j \rangle \langle b_j \, ^{\cdot} v, b_i \rangle a_i 
\vspace{3mm}\\
& \ds \quad - \langle u, b_j \rangle \langle a_j \, ^{\cdot} v, a_i \rangle b_i 
+ \langle u, a_j \rangle \langle b_j \, ^{\cdot}  v, a_i \rangle b_i ) 
\vspace{3mm}\\
& \ds \quad = \sum_{ij} (\langle u, b_j\rangle \langle v, [b_i, a_j] \rangle a_i 
- \langle u, a_j \rangle \langle v, [b_i, b_i] \rangle a_i 
\vspace{3mm}\\
& \ds \quad - \langle u, b_j \rangle \langle u, [a_i, a_j] \rangle b_i + 
\langle u, a_j \rangle \langle v, [a_i, b_j] \rangle b_i ) 
\vspace{3mm}\\
& \quad \ds {\mathop{=}\limits^{\mbox{\scriptsize [by (1.29b)]}}} \ 
 -(\langle  u, \rangle \otimes \langle v, \rangle \otimes {\mathbf 1}) ([r^{12}, r^{23}]); 
\ea
\ee
Interchanging $u$ and $v$ in the above calculation, we get 
\renewcommand{\theequation}{\arabic{section}.\arabic{equation}{\rm b}}
\setcounter{equation}{21}
\be
\ba{ll}
2) &  \ds - {\mathcal O} ({\mathcal O} (v)^. u) = 
 \sum_{ij} (\langle u, [a_j, b_i)] \rangle \langle v, b_j \rangle a_i + \langle u, [b_i, b_j] 
\rangle \langle v, a_j \rangle a_i 
\vspace{3mm}\\
& \ds \quad + \langle u, [a_i, a_j] \rangle \langle v, b_j \rangle b_i 
- \langle u, [a_i, b_j] \rangle \langle v, a_j \rangle b_i) 
\vspace{3mm}\\
& \ds \quad {\mathop{=}\limits^{\mbox{\scriptsize [by (1.29a)]}}} \
 - (\langle u, \rangle \otimes \langle v, \rangle \otimes {\mathbf 1}) ( [r^{12}, r^{13}]); 
\ea
\ee
\renewcommand{\theequation}{\arabic{section}.\arabic{equation}{\rm c}}
\setcounter{equation}{21}
\be
\ba{ll}
3) & \ds  [{\mathcal O} (v), {\mathcal O} (u)] = \sum_{ij} [\langle v, b_i \rangle a_i 
- \langle v, a_i \rangle b_i, \langle u, b_j \rangle a_j - 
\langle u, a_j \rangle b_j] 
\vspace{3mm}\\
& \ds \quad = \sum_{ij} (\langle u, b_j\rangle \langle v, b_i\rangle  [a_i, a_j] 
- \langle u, a_j \rangle \langle v, b_i\rangle [a_i, b_j] 
\vspace{3mm}\\
& \ds \quad - \langle u, b_j \rangle \langle v, a_i\rangle [b_i, a_j ] + 
\langle u, a_j\rangle \langle v, a_i\rangle [b_i, b_j] 
\vspace{3mm}\\
& \ds \quad {\mathop{=}\limits^{\mbox{\scriptsize [by (1.29c)]}}} \
 - (\langle u, \rangle \otimes \langle v, \rangle \otimes {\mathbf 1}) ([r^{13},  r^{23}]).
\ea
\ee
Altogether, we thus f\/ind
\renewcommand{\theequation}{\arabic{section}.\arabic{equation}}
\setcounter{equation}{22}
\be
\langle w, [{\mathcal O} (u), {\mathcal O} (v)] - {\mathcal O} ({\mathcal O} (u)^. v 
- {\mathcal O} (v)^. u) \rangle = 
\langle u \otimes v \otimes w,  c(r)\rangle , \quad
 \forall \; u, v, w \in {\mathcal G}^*. 
\ee

\setcounter{equation}{0}

\section{Dif\/ferential Lie Algebras say:  2-cocycles,\\
\hspace*{12mm}-- Si, {\mathversion{bold}$r$}-matrices, -- No,
{\mathversion{bold}${\mathcal O}$}-operators are welcome}

Let $R$ be a commutative ring or algebra, and $\partial_1, \ldots, \partial_m,:  R \rightarrow 
R$ be $m$ commuting derivations.  A Lie algebra over $R$ is $R^N$, some $N \in{\mathbf{N}}$, 
with a  skew multiplication [~,~] $R^N \times R^N \rightarrow R^N$ given by bilinear {\it differential
operators}, such that
\be
[[x, y], z] + \mbox{c.p.} = 0 , \qquad \forall \; x, y, z \in \tilde R^N, 
\ee
where $\tilde R \supset R$ is {\it arbitrary} dif\/ferential extension of $R$.  (This means that 
the skewsymmetry of [~,~] and Jacobi identity are the properties solely of dif\/ferential operators 
performing the multiplication [~,~] and are not dependent upon the quirks of $R$ itself).

In this Section we consider the very simplest  case $m = 1$.  Denote $\partial_1$ by $\partial$.  
First, let ${\mathcal D}_1$ be the Lie algebra of vector f\/ields on the line:  ${\mathcal D}_1 = R$ 
and 
\be
[X, Y] = XY^\prime - X^\prime Y, \qquad  \forall \; X, Y \in {\mathcal D}_1, 
\ee
where $( \cdot )^\prime: = \partial ( \cdot )$.

Consider the following bilinear form on ${\mathcal D}_1$:
\be
\omega (X, Y) = \partial^3 (X) Y. 
\ee
This form is skewsymmetric, in {\it differential  sense}, for
\be
\omega (X, Y) \sim - \omega (Y, X), \qquad \forall \; X, Y \in {\mathcal D}_1, 
\ee
where $a \sim b$ means that $(a - b) \in \mbox{Im} \; \partial.$  Also, 
\be
\omega (X, [Y, Z]) + \mbox{c.p.} \sim 0, \qquad \forall \; X, Y, Z \in {\mathcal D}_1, 
\ee
so that $\omega$ is a (generalized) 2-cocycle.  (All the necessary details of the theory 
can be found in~[6].)  But if we try to represent this 2-cocycle $\omega$ in the 
${\mathcal O}$-form (2.3),
\be
\omega (X, Y) = \langle {\mathcal O}^{-1} (X), Y \rangle, 
\ee
we f\/ind that ${\mathcal O} = \partial^{-3}$: in other words, ${\mathcal O}$ doesn't exist, 
and neither does~$r$.  We circumvent this particular obstacle as follows.

Denote by ${\mathcal G} (\mu), \ \mu = \mbox{const} \in {\mathcal F}: = \mbox{Ker} \; \partial|_R$, 
the following Lie algebra structure 
on $R + R$:
\be
\left[ \left(\matrix{X \cr f \cr} \right),  \left(\matrix{Y \cr g \cr} \right) \right] = 
\left(\matrix{XY^\prime - X^\prime Y \cr (Xg - Yf + \mu (X^\prime Y^{\prime \prime} - X^{\prime 
\prime} Y^\prime )) ^\prime \cr} \right). 
\ee
We still have the 2-cocycle $\omega$ (3.3) on ${\mathcal G} (\mu)$:
\be
\omega \left( \left(\matrix{X \cr f \cr} \right), \left(\matrix{Y \cr g \cr} \right) \right) = 
\partial^3 (X) Y, 
\ee
and it is still degenerate.  However, ${\mathcal G}(\mu)$ also possess a {\it nondegenerate}  
symplectic 2-cocycle.
\be
\Omega \left(  {X \choose f} , {Y \choose g} \right) =  Xg - Yf. 
\ee
Indeed,
\[
\ba{l}
\ds \Omega \left( {X \choose f},  \left[ {Y \choose g}, {Z \choose h} \right]  \right) + 
\mbox{c.p.}
\vspace{3mm}\\
\ds \qquad 
= \Omega \left( {X \choose f}, \left( {YZ^\prime - Y^\prime Z \choose (Yh - Zg + \mu (Y^\prime 
Z^{\prime \prime} - Y^{\prime \prime} Z^\prime))^\prime } \right) \right)+ \mbox{c.p.}
\vspace{3mm}\\
\ds \qquad = \left\{ - (YZ^\prime - Y^\prime Z) f + X (Y^\prime h + Yh^\prime - Z^\prime g 
- Zg^\prime + \mu Y^\prime Z^{\prime \prime \prime}  - \mu Y^{\prime \prime \prime} Z^\prime) 
\right\} + \mbox{c.p.}
\vspace{3mm}\\
\ds \qquad = \mu X (Y^\prime Z^{\prime \prime \prime} - Y^{\prime \prime \prime} Z^\prime)
 + \mbox{c.p.}
\sim - \mu X^\prime (Y^\prime Z^{\prime \prime \prime} - Y^{\prime \prime \prime} Z^\prime) + 
\mbox{c.p.} = 0. 
\ea
\]
Taking the sum $\epsilon \omega + \Omega$ as the new 2-cocycle on ${\mathcal G}(\mu)$, 
$\epsilon = \mbox{const}$, we get
\be
\ba{l}
\ds (\epsilon \omega + \Omega) \left( {X \choose f}, {Y \choose g} \right) 
\vspace{3mm}\\
\ds \qquad = \epsilon \partial^3 (X) Y + (Xg - Yf) = 
 \langle \left(\matrix{\epsilon \partial^3 & -1 \cr 1 & 0 \cr} \right) \left(\matrix{ X \cr f \cr} 
\right), \left(\matrix{Y \cr g \cr} \right) \rangle, 
\ea
\ee
so that
\be
{\mathcal O}^{-1} = \left(\matrix{\epsilon \partial^3 & -1 \cr 1 & 0 \cr} \right), 
\ee
and thus
\be
{\mathcal O} = \left(\matrix{0 & 1 \cr -1 & \epsilon \partial^3 \cr} \right). 
\ee
 
Since ${\mathcal O}$ is a {\it differential  operator}, it corresponds to no element of 
$\wedge^2 {\mathcal G} (\mu)$; it is only in f\/inite dimensions that we can identify 
${\mathcal G}$ with $\mbox{Hom}\, ({\mathcal G}^*, R)$ and 
${\mathcal G}^{\otimes 2}$ with $\mbox{Hom}\,({\mathcal G}^{* \otimes 2}, R)$.  
The conclusion is inescapable:  the 
proper notion of the classical $r$-matrix is a skewsymmetric ${\mathcal O}$-operator 
${\mathcal G}^* \rightarrow 
{\mathcal G}$ satisfying the classical ${\mathcal O}$-def\/ining equation~(2.10).

In f\/inite dimensions, there exists a dif\/ferent version of the notion of classical $r$-matrix, 
due to Semeynov-Tyan-Shansky~[13].  It is already in operator  form, acting as $r: {\mathcal G} 
\rightarrow {\mathcal G}$; but it requires ${\mathcal G}$ to have an invariant 
nondegenerate scalar product, a 
condition rarely encountered in dif\/ferential situations.  [E.g., the Lie algebra 
${\mathcal G} (\mu)$~(3.7) has no invariant scalar product, no matter what $\mu$ is.) 

We conclude this Section by calculating the Lie algebra structure on ${\mathcal G}(\mu)^*$ 
induced by the ${\mathcal O}$-operator~(3.12) via formula~(2.11).

Denote typical elements in ${\mathcal G} (\mu)^*$ as ${{\ds{u \choose p}, {v \choose q}}},$ 
with the pairing
\be
\langle {u \choose p},  {X \choose f} \rangle = u X + p f. 
\ee
Let us f\/irst obtain the formula for the coadjoint action of ${\mathcal G} (\mu) $ 
on ${\mathcal G} (\mu)^*$:
\[
\ba{l}
\ds \langle {X \choose f} ^. {u \choose p},  {Y \choose g} \rangle  \sim - \langle {u \choose p}, 
\left[ {X \choose f}, {Y \choose g} \right] \rangle
\vspace{3mm}\\
 \ds \qquad {\mathop{\sim }\limits^{\mbox{\scriptsize [by (3.7)]}}} \ 
 - u (XY^\prime - X^\prime Y) + p^\prime (Xg - Yf + \mu X^\prime Y^{\prime \prime} - 
\mu X^{\prime \prime} Y^\prime) 
\vspace{3mm}\\
\ds \qquad \sim ((u X)^\prime + u X^\prime) Y - fp^\prime Y + X p^\prime g + ((\mu X^\prime p^\prime)
^{\prime \prime} + (\mu X^{\prime \prime} p^\prime)^\prime) Y. 
\ea
\]
 
Thus,
\be
{X \choose f}^.  {u \choose p} = \left(\matrix{X \partial + 2 X^\prime & 
(\mu \partial (X^\prime \partial + 2 X^{\prime \prime}) - f) \partial \cr
0 & X \partial \cr} \right) {u \choose p}. 
\ee

Hence,
\be
\ba{l}
\ds {\mathcal O} {u \choose p}^. {v \choose q} = {p \choose -u + \epsilon p^{\prime \prime \prime}} ^.
{v \choose q} 
\vspace{3mm}\\
\ds \qquad = {(pv^\prime + 2 p^\prime v) + \mu ((p^\prime q^\prime)^{\prime \prime} + (p^{\prime \prime} 
q^\prime)^\prime) - (- u + \epsilon p^{\prime \prime \prime}) q^\prime \choose pq^\prime }.
\ea
\ee
Therefore,
\be
\ba{l}
\ds \left[ {u \choose p}_, {v \choose q} \right] = {\mathcal O} {u \choose p}^. {v \choose q}  
- {\mathcal O} {v \choose q}^. {u \choose p} 
\vspace{3mm}\\
\ds \qquad ={(pv - qu + (\epsilon - \mu) (p^\prime q^{\prime \prime} - p^{\prime \prime} q^{\prime}))
^\prime \choose pq^\prime - p^\prime q} . 
\ea
\ee
We see that
\be
{\mathcal G} (\mu)^* \approx {\mathcal G} (\epsilon - \mu). 
\ee
Since ${\mathcal G}(0) = {\mathcal D}_1 \ltimes V_1$ is certainly a Lie algebra, 
${\mathcal G} (\mu)$ is so aposteriori:
\be
{\mathcal G} (\epsilon) \approx {\mathcal G} (0)^*. 
\ee
The reader who hasn't  bothered to check the Jacobi identity for the Lie bracket~(3.7) on 
${\mathcal G}(\mu)$ may now feel smug about it.  The reader who didn't blink an 
eye when the symplectic 2-cocycle~(3.9) 
was sprung out on him as a deus ex machina without an explanation, as this 
is how modern mathematics is supposed to operate, will be disappointed to f\/ind a general 
construction of symplectic $r$-matrices in Appendix~A2.  Sorry about that.

\setcounter{equation}{0}

\section{ {\mathversion{bold}${\mathcal O}$}-natural property of the 
{\mathversion{bold}${\mathcal O}$}-operators}

Let $\varphi: {\mathcal G} \rightarrow {\mathcal H}$ be a homomorphism of Lie algebras.  
If everything is f\/inite-dimensional and $r \in \wedge^2 {\mathcal G}$ then 
$\varphi (r) \in \wedge^2 {\mathcal H}$, 
and 
\be
c (\varphi (r)) = \varphi (c(r)); 
\ee
thus, if $r$ is a classical $r$-matrix then so is $\varphi (r)$.

Consider now the general case.  
Let ${\mathcal O} = {\mathcal O}_{\mathcal G}: {\mathcal G}^* \rightarrow {\mathcal G}$ be an 
${\mathcal O}$-operator.  Recall that this means that 
\be
{\mathcal O} ({\mathcal O} (u)^. v - {\mathcal O} (v)^. u) = [{\mathcal O} (u), {\mathcal O}(v)],
\qquad  \forall \; u, v \in {\mathcal G}^*, 
\ee
and that ${\mathcal O}$ is skewsymmetric:
\be
\langle u, {\mathcal O} (v) \rangle
 + \langle v, {\mathcal O} (u) \rangle  \sim 0, \qquad  \forall \;u, v \in {\mathcal G}^*. 
\ee
Since
\be
\langle u, {\mathcal O} (v) \rangle + 
\langle v, {\mathcal O} (u)\rangle  \sim 
\pbf{u}^t ({\mathcal O} + {\mathcal O}^\dagger) (\pbf{v}), 
\ee
${\mathcal O}$ is skewsymmetric if\/f
\be
{\mathcal O}^\dagger = - {\mathcal O}, 
\ee
where ${\mathcal O}^\dagger$ is the operator adjoint to ${\mathcal O}$, and $\pbf{u}$ and 
$\pbf{v}$ are treated as column-vectors (see~[6].)

Now let $\varphi: {\mathcal G} \rightarrow {\mathcal H}$ be a homomorphism of  Lie algebras.  
It induces the dual map
$\varphi^* : {\mathcal H}^* \rightarrow {\mathcal G}^*$.  Since
\be
\langle \varphi^* (\bar u), x\rangle = \langle {\bar u}, \varphi (x)\rangle, 
\qquad  \forall\; {\bar u} \in {\mathcal H}^*, \quad \forall \; x \in {\mathcal G}, 
\ee
and
\be
\langle {\bar u}, \varphi (x) \rangle  = {\bar {\pbf{u}}}^t \varphi (x) \sim (\varphi^\dagger 
({\bar {\pbf{u}}}))^t x, 
\ee
we see that
\be
\varphi^* = \varphi^\dagger. 
\ee

\noindent
{\bf Proposition 4.9.} {\it  Set 
\setcounter{equation}{9}
\be
{\mathcal O}_{\mathcal H} = 
\varphi {\mathcal O}_{\mathcal G} \varphi^*: {\mathcal H}^* \rightarrow {\mathcal H}.
\ee
Then ${\mathcal O}_{\mathcal H}$ is an ${\mathcal O}$-operator.}

\medskip

\noindent{\bf Proof.} 
 First, let's check that ${\mathcal O}_{\mathcal H}$ is skewsymmetric.  We have:
\be
({\mathcal O}_{\mathcal H})^\dagger = 
(\varphi {\mathcal O}_{\mathcal G} \varphi^\dagger)^\dagger = \varphi 
( {\mathcal O}_{\mathcal G})^\dagger \varphi^\dagger =
 - \varphi {\mathcal O}_{\mathcal G} \varphi^\dagger = - {\mathcal O}_{\mathcal H}. 
\ee
Next, for any $\bar u, \bar v \in {\mathcal H}^*$, we have to verify that
\be
{\mathcal O}_{{\mathcal H}} ({\mathcal O}_{{\mathcal H}} (\bar u)^. \bar v 
- {\mathcal O}_{{\mathcal H}} (\bar v)^. \bar u) = [{\mathcal O}_{{\mathcal H}} (\bar u), 
{\mathcal O}_{{\mathcal H}} (\bar v) ]. 
\ee
For the LHS of the expression (4.12) we get
\renewcommand{\theequation}{\arabic{section}.\arabic{equation}{\rm L}}
\setcounter{equation}{12}
\be
\varphi {\mathcal O}_{\mathcal G} \varphi^\dagger (\varphi {\mathcal O}_{\mathcal G}^\dagger 
(\bar u)^. \bar v - \varphi {\mathcal O}_{\mathcal G} \varphi^\dagger (\bar v)^. \bar u ), 
\ee
and for the RHS of the expression (4.12) we obtain
\renewcommand{\theequation}{\arabic{section}.\arabic{equation}{\rm R}}
\setcounter{equation}{12}
\be
\ba{l}
[\varphi {\mathcal O}_{\mathcal G} \varphi^\dagger (\bar u), 
\varphi {\mathcal O}_{\mathcal G} \varphi^\dagger (\bar v)] \
{\mathop{=}\limits^{\mbox{\scriptsize [since $\varphi$ is  a  homomorphism]}}} \
 \varphi [ {\mathcal O}_{\mathcal G} \varphi^\dagger (\bar u), {\mathcal O}_{\mathcal G} 
\varphi^\dagger (\bar v)] \ 
\vspace{2mm}\\
\ds \qquad {\mathop{=}\limits^{\mbox{\scriptsize [since  ${\mathcal O}_{\mathcal G}$ 
is an  ${\mathcal O}$-operator]}}} \
 \varphi {\mathcal O}_{\mathcal G} ({\mathcal O}_{\mathcal G} 
\varphi^\dagger (\bar u)^. \varphi^\dagger (\bar v) - {\mathcal O}_{\mathcal G} 
\varphi^\dagger (\bar v)^. \varphi^\dagger (\bar u)), 
\ea
\ee
and by formula (4.15) below the expressions (4.13$L,R$) are equal. \hfill \rule{3mm}{3mm}

\medskip

\renewcommand{\theequation}{\arabic{section}.\arabic{equation}}
\setcounter{equation}{14}

\noindent
{\bf Lemma 4.14.} 
\be
\varphi^\dagger (\varphi (X)^. \bar v) = X^. \varphi^\dagger (\bar v), \qquad  \forall\;
X \in {\mathcal G}, \quad  \forall \; \bar v \in {\mathcal H}^*. 
\ee

\noindent
{\bf Proof.}  Formula (4.15) is an equality in ${\mathcal G}^*$.  Any such equality, 
$(\cdot) = (\cdot \cdot)$, is equivalent to the relation 
\be
\langle (\cdot), Y \rangle \sim  \langle (\cdot \cdot), Y \rangle, \qquad 
\forall \; Y \in {\mathcal G}. 
\ee
So, 
\[
\ba{l}
\ds \langle \varphi^\dagger (\varphi (X)^. \bar v), Y \rangle  \sim 
\langle \varphi (X)^. \bar v, \varphi (Y)\rangle  \sim - \langle \bar v, [\varphi (X), \varphi (Y)] \rangle
\vspace{2mm}\\
\ds {\mathop{=}\limits^{\mbox{\scriptsize [since $\varphi$ is a homomorphism]}}}
\ - \langle \bar v, \varphi ([X, Y]) \rangle  \sim - \langle \varphi^\dagger (\bar v), [X, Y] \rangle\sim 
 \langle X^.  \varphi^\dagger (\bar v), Y\rangle. \hspace{28.8pt}  \mbox{\rule{3mm}{3mm}}
\ea
\] 

In \S~6 we establish that a quadratic Poisson bracket on ${\mathcal G}^*$ canonically 
associated to 
every ${\mathcal O}$-operator ${\mathcal O}_{\mathcal G}$, also has the natural property.  

\medskip

\noindent
{\bf Remark 4.17.} The map $\varphi^*: {\mathcal H}^* \rightarrow {\mathcal G}^*$ 
is a homomorphism of Lie algebras.

\medskip

\noindent
{\bf Proof.}  Take any $\bar u, \bar v \in {\mathcal H}^*$.  We have to show that 
\setcounter{equation}{17}
\be
\varphi^\dagger ([\bar u, \bar v)] = [\varphi^\dagger (\bar u), \varphi^\dagger (\bar v)], 
\ee
which is
\[
\varphi^\dagger ({\mathcal O}_{{\mathcal H}} (\bar u)^. \bar v 
- {\mathcal O}_{{\mathcal H}}  (\bar v)^. \bar u) = 
{\mathcal O}_{\mathcal G} \varphi^\dagger (\bar u)^. \varphi^\dagger (\bar v) 
- {\mathcal O}_{\mathcal G} \varphi^\dagger (\bar v)^. \varphi^\dagger (\bar u),
\]
which further is
\[
\varphi^\dagger (\varphi {\mathcal O}_{\mathcal G} \varphi^\dagger (\bar v)^. \bar v 
- \varphi {\mathcal O}_{\mathcal G} \varphi^\dagger (\bar v)^. \bar u) = 
{\mathcal O}_{\mathcal G} \varphi^\dagger (\bar u)^. \varphi^\dagger (\bar v)
 - {\mathcal O}_{\mathcal G} \varphi^\dagger (\bar v)^. \varphi^\dagger (\bar u), 
\]
and this is true by formula (4.15). \hfill \rule{3mm}{3mm}

\medskip

\noindent
{\bf Remark 4.19.}  In the generality we are working, many f\/inite-dimensional notions 
disappear.  For example, Hamilton--Lie groups and their Hamiltonian actions, -- though 
{\it infinitestimal} versions of those may survive, see \S~5,~6.   Other things disappear 
altogether, such as representation of commutator on ${\mathcal G}^*$ by the cocommutator 
${\mathcal G} \rightarrow {\mathcal G} \otimes {\mathcal G}$.  Every time one needs 
to use the f\/inite-dimensional isomorphism $\mbox{End}\, (V) \approx 
V^* \otimes V$, one gets into all sorts of trouble with the authorities, since 
$\mbox{Dif\/f}\,(V)$ is inf\/inite-dimensional no matter what dimension of $V$ is.

\setcounter{equation}{0}

\section{Linear Poisson brackets on dual spaces to Lie algebras}

Before we tackle {\it quadratic} Poisson brackets on ${\mathcal G}^*$, 
it's instructive to review the 
{\it linear} Poisson brackets; this way we can introduce basic def\/initions and themes in 
more familiar surroundings.

So, let ${\mathcal G} = R^N$ be a dif\/ferential Lie algebra.  
(Or dif\/ferential-dif\/ference one; it's almost 
the same, as far as the theory goes, so I prefer not to clutter the presentation with 
indices corresponding to discrete degrees of freedom.  See Remark~5.50 below for more details.)  
The dual space ${\mathcal G}^*$ is also $R^N$.  The dif\/ferential ring 
$C = C_u = R[u_i^{(\sigma)}]$, $i = 1,\ldots, N$, $\sigma \in {\mathbf{Z}}^m_+,$ 
is what used to be the ring $\mbox{Fun}\,({\mathcal G}^*)$ of smooth 
functions on ${\mathcal G}^*$ in f\/inite dimensions.

On the ring $C_u$ we have the Poisson bracket
\be
\{H, F\} = X_H (F) \sim {\delta F \over \delta \pbf{u}^t} B \left( {\delta H \over \delta \pbf{u}} 
\right), 
\ee
where the Hamiltonian matrix $B$, {\it linear} in $u$, is extracted from the following 
def\/ining relation:
\be
\{ \pbf{u} ^t \pbf{X}, \pbf{u} ^t \pbf{Y}\} \sim \pbf{u} ^t [\pbf{X}, \pbf{Y}], \qquad
 \forall \; \pbf{X}, \pbf{Y} \in {\mathcal G}. 
\ee
This is the dif\/ferential version of the more familiar form
\be
\{\langle u, X \rangle, \langle u, Y\rangle \} \sim  \langle u, [X, Y] \rangle, 
\qquad \forall \; X, Y \in {\mathcal G}. 
\ee
Since
\renewcommand{\theequation}{\arabic{section}.\arabic{equation}{\rm a}}
\setcounter{equation}{3}
\be
\{\pbf{u}^t \pbf{X}, \pbf{u}^t \pbf{Y}\} \sim \pbf{Y}^t B (\pbf{X}),
\ee
and 
\renewcommand{\theequation}{\arabic{section}.\arabic{equation}{\rm b}}
\setcounter{equation}{3}
\be
\pbf{u}^t [\pbf{X}, \pbf{Y}] = \langle \pbf{u}, [ \pbf{X}, \pbf{Y}] \rangle  \sim - 
\langle \pbf{X}^. \pbf{u}, \pbf{Y} \rangle, 
\ee
we see that
\renewcommand{\theequation}{\arabic{section}.\arabic{equation}}
\setcounter{equation}{4}
\be
B (\pbf{X}) = - \pbf{X}^. \pbf{u}.
\ee
Thus, the linear Poisson bracket on ${\mathcal G}^*$ has the form
\renewcommand{\theequation}{\arabic{section}.\arabic{equation}{\rm a}}
\setcounter{equation}{5}
\be
\{ H, F\} \sim - \langle {\delta H \over \delta \pbf{u}}^. \pbf{u}, {\delta F \over \delta \pbf{u}}\rangle
\ee
\renewcommand{\theequation}{\arabic{section}.\arabic{equation}{\rm b}}
\setcounter{equation}{5}
\be
\qquad \sim \langle \pbf{u}, \left[{\delta H \over \delta \pbf{u}}, 
{\delta F \over \delta \pbf{u}} \right]. 
\ee
The {\it Casimirs} of a Poisson bracket are those Hamiltonians $H$ for which the vector 
f\/ield $X_H = \{H, . \} $ is identically zero.  (In f\/inite dimensions, the common level 
surfaces of Casimirs are symplectic leaves.)  From formula (5.6a) 
we see that the Casimirs on 
${\mathcal G}^*$, also called (for a reason) coadjoint invariants, are the solutions of the equation
\renewcommand{\theequation}{\arabic{section}.\arabic{equation}}
\setcounter{equation}{6}
\be
{\delta H \over \delta \pbf{u}}^. \pbf{u} = {\mathbf 0}.
\ee
Equivalently,
\be
\langle \pbf{X}^. \pbf{u}, {\delta H \over \delta \pbf{u}} \rangle  \sim 0, \qquad
 \forall \; X \in {\mathcal G}. 
\ee

\noindent
{\bf Example 5.9.}  Let ${\mathcal G}$ be the Lie lagebra ${\mathcal D}_1$ of \S~3:
\setcounter{equation}{9}
\be
[X, Y] = X \partial (Y) - \partial (X) Y, \qquad \forall\; X, Y \in {\mathcal D}_1. 
\ee
Then
\be
u[X, Y] = u (XY^\prime - X^\prime Y) \sim - Y (u \partial + \partial u) (X). 
\ee
Thus, 
\be
B = B ({\mathcal D}_1) = - (u \partial + \partial u) 
\ee
\be
\qquad = - 2 {\sqrt u} \partial {\sqrt u}. 
\ee
Therefore,
\be
H \in \mbox{Ker}\, (B) \ \Leftrightarrow \ {\delta H \over \delta u} =  \mbox{const}/ \sqrt{u }
\ \Leftrightarrow \ H = \mbox{const} \sqrt{u}. 
\ee
Equivalently, from formula (3.14) we see that
\be
X^. u = (X \partial + 2 X^\prime) (u) = Xu^\prime + 2X^\prime u = {1\over X} 
(X^2 u)^\prime, 
\ee
so that $H$ is Casimir if\/f
\be
\left({\delta H \over \delta u} \right)^2 u = \mbox{const} \ \Leftrightarrow \ H = 
\mbox{const}\sqrt u. 
\ee
We see that in this case $H$ belongs not to the ring $C_u$ itself but to its algebraic 
extension.

\medskip

\noindent
{\bf Remark 5.17.}  In f\/inite dimensions, the linear Poisson bracket (5.6) was 
discovered by Lie and rediscovered by everyone else.

\medskip

Let us check that the linear Poisson bracket (5.6) is {\it natural}.  
Let $\varphi: {\mathcal G} \rightarrow 
{\mathcal H}$ be a homomorphism of Lie algebras.  Let $C_q = R[q_j^{(\sigma)}]$ 
be the dif\/ferential ring of functions on~${\mathcal H}^*$, 
$j = 1, \ldots, \mbox{dim}\,({\mathcal H}).$  Let $\Phi: C_u \rightarrow 
C_q$ be the dif\/ferential homomorphism dual to the map of spaces
 $\varphi^*: {\mathcal H}^* \rightarrow {\mathcal G}^*$.  
To calculate $\Phi$, it's enough to notice that $C_u$ and $C_q$ are (dif\/ferentially) 
generated by linear functions on ${\mathcal G}^*$ and ${\mathcal H}^*$ respectively:  
\setcounter{equation}{17}
\be
\Phi:  \langle \ , X \rangle  \longmapsto  \langle \ , \varphi (X) \rangle, \qquad
 \forall \; X \in {\mathcal G}. 
\ee
Thus,
\be
\Phi ( \langle u, X \rangle ) = \Phi (\pbf{u})^t \pbf{X} = 
\langle \pbf{q},  \varphi (\pbf{X}) \rangle = \pbf{q}^t \varphi (\pbf{X}) 
\sim \varphi^\dagger (\pbf{q})^t \pbf{X}, 
\ee
so that 
\be
\Phi (\pbf{u}) = \varphi^\dagger (\pbf{q}).
\ee
Denote by $B_{\mathcal G} = B({\mathcal G})$ the linear Hamiltonian matrix associated 
to the Lie algebra ${\mathcal G}$ by 
formula (5.5).  The property of the matrix $B_{\mathcal G}$ being natural means that
\be
\Phi (\{H, F\}_{{\mathcal G}^{*}}) = \{ \Phi (H), \Phi (F) \}_{{\mathcal H}{^*}}, \qquad 
\forall \; H, F \in C_{{\mathcal G}^{*}} = C_u.
\ee
By the well-known criterion (see, e.g.~[6] p.~54), 
a map $\Phi:  C_1 \rightarrow C_2$ is 
Hamiltonian between the Hamiltonian matrices $B_1$ and $B_2$ over rings $C_2$ and $C_2$ 
respectively, if\/f
\be
\Phi (B_1) = D ({\pbf{\Phi}}) B_2 D ({\pbf{\Phi}})^\dagger, 
\ee
where $D$ stands for the Fr\`echet derivative and 
\be
{\pbf{\Phi}} = \Phi (\pbf{q}^1),
\ee
$\pbf{q}^1$ and $\pbf{q}^2$ being the column-vectors of generators of the rings $C_1$ and $C_2$ 
respectively.  For the matrices $B_1 = B_{\mathcal G} $ and $B_2 = B_{{\mathcal H}}$, 
a proof of the identity~(5.22) can be found in~[6] p.~66.  
Instead of repeating this type of proof, -- which becomes 
very cumbersome for the quadratic Poisson bracket on ${\mathcal G}^* $ def\/ined in 
\S~6, -- I will 
reformulate the criterion (5.22) into a very useful form:

\medskip

\noindent
{\bf Proposition 5.24.} {\it  A map $\Phi: C_1 \rightarrow C_2$ is Hamiltonian iff
\setcounter{equation}{24}
\be
\Phi (\{H, F\}_1) \sim \{ \Phi (H), \Phi (F) \}_2, \qquad
 \forall \; H, F \quad  \mbox{linear in} \ \pbf{q}^1. 
\ee}

\noindent
{\bf Proof.} We shall show that if $H$ and $F$ are {\it arbitrary  linear} in 
$\pbf{q}^1$ then formula (5.25) implies formula (5.22).  Let 
\[
H = \pbf{q}^{1t} \pbf{X}, \qquad F = \pbf{q}^{1t} \pbf{Y}, \qquad 
\pbf{X}, \pbf{Y} \in \tilde R^{N_{1}},
\]
so that
\be
\{H, F\}_1 \sim \pbf{Y}^t B_1 (\pbf{X}), \qquad \Phi (\{H, F\}_1) \sim \pbf{Y}^t \Phi 
(B_1) (\pbf{X}), 
\ee
\be
\{ \Phi (H), \Phi (F) \}_2 = \{ \pbf{\Phi}^t \pbf{X}, \pbf{\Phi}^t \pbf{Y}\}_2 \sim \left( {\delta \over 
\delta \pbf{q}^2} (\pbf{\Phi}^t \pbf{Y})\right)^t B_2 {\delta \over \delta \pbf{q}^2} 
(\pbf{\Phi}^t \pbf{X}). 
\ee
Now, since $\pbf{X}$ and $\pbf{Y}$ are $q^{(.)}$-independent, 
\be
{\delta \over \delta q_s} (\pbf{\Phi}^t \pbf{X}) = D_{q_{s}} (\pbf{\Phi})^\dagger  (\pbf{X}) 
= \sum_j D_{q_{s}} (\Phi_j)^\dagger (X_j), 
\ee
so that
\be
\ba{l}
\ds \left({\delta \over \delta \pbf{q}^2} (\pbf{\Phi}^t \pbf{Y})\right)^t B_2 
{\delta \over \delta \pbf{q}^2} (\pbf{\Phi}^t 
\pbf{X}) = \sum_{ij} [D_i (\pbf{\Phi})^\dagger (\pbf{Y})]^t (B_2)_{ij} D_j (\pbf{\Phi})^\dagger 
(\pbf{X}) 
\vspace{3mm}\\
\ds \qquad \sim {\sum_{ij}} \pbf{Y}^t D_i (\pbf{\Phi}) (B_2)_{ij} 
D_j (\pbf{\Phi})^\dagger (\pbf{X}) = \pbf{Y}^t D 
(\pbf{\Phi}) B_2 D (\pbf{\Phi})^\dagger (\pbf{X}).
\ea
\ee
Comparing formulae (5.26) and (5.29) and remembering that $\pbf{X}$ and 
$\pbf{Y}$ are arbitrary, we arrive at the Hamiltonian criterion (5.22). \hfill \rule{3mm}{3mm}

\medskip

\noindent
{\bf Remark 5.30.} A little more ef\/fort will show that the equality modulo 
$\sum\limits^m_{\ell=1} \mbox{Im}\,  \partial_\ell$ sign~$\sim$ 
in formula (5.25) can be replaced by the exact 
equality sign = (see~[6] p.~53).  We won't need this more precise form in what follows.

\medskip

Assume for a moment that we are thrown back in time into f\/inite dimensions.  Let $G$ be a Lie 
group whose Lie algebra is ${\mathcal G}$.  Then $G$ acts on ${\mathcal G}^*$ 
by the coadjoint representation, 
and this action preserves the linear Poisson bracket on ${\mathcal G}^*$.  
This means that the map 
$\mbox{Ad}^*: G \times {\mathcal G}^* \rightarrow {\mathcal G}^*$ is Poisson, 
with the Poisson bracket on $G \times 
{\mathcal G}^*$ being the product of Poisson brackets on $G$ 
and ${\mathcal G}^*$ and Poisson bracket and $G$ 
being zero.  This is a particular case of the following more general set-up.  Let $G$ be a 
Hamilton--Lie group.  (This is the original name given to the subject by its inventor, 
V.G.~Drinfel'd~[2]; subsequent commentators have changed the original name into 
``Poisson--Lie'' groups).  Let $M$ be a Poisson manifold.  
Suppose $G$ acts from the left on $M$ in
such a way that the action map $G \times M \rightarrow M$ is a Hamiltonian (=~Poisson) map, 
with the Poisson structure on $G \times M$ being of product type.  Then {\it infinitesimal} 
criterion for this action to be Hamiltonian is~([14])
\setcounter{equation}{30}
\be
\ba{l}
\ds X^\wedge (\{H, F\}) - \{ X^\wedge (H), F\} - \{H, X^\wedge (F) \} = \langle [\theta _H, \theta_F], 
X \rangle,
\vspace{2mm}\\
\ds \forall \; H, F \in \mbox{Fun}\, (M), \qquad  \forall \; X \in {\mathcal G} = \mbox{Lie}\, (G).
\ea
\ee
Here $X^\wedge$ is the vector f\/ield $\ds {d \over dt} \exp (tX)^*\bigg|_{t=0}$ on $M$ generated by 
$X \in {\mathcal G}$, $\{ \ ,  \}$ is the Poisson bracket on $M$, $ \theta_H:  M \rightarrow 
{\mathcal G}^*$, for a given function $H$ on $M$, is the map def\/ined by the rule
\be
\theta_H (x) = d_g H (gx)|_{g = e},  \qquad \forall\; x \in M, 
\ee
and $[\theta_H, \theta_F]$ is the commutator in ${\mathcal G}^*$ induced by the dif\/ferential at the 
identity in $G$ of the multiplicative Poisson bracket on $G$.  (For a proof, see, e.g.~[4] p.~45).

We aim to reformulate the inf\/initesimal criterion of Hamiltonian action (5.31) into a 
{\it definition} usable for the functional case where there are no Lie groups present 
anymore, only Lie algebras.  Let ${\mathcal G}$ 
be such Lie algebra, and let $^\wedge: {\mathcal G} \mapsto D^{ev}(C)$ 
be an antirepresentation of ${\mathcal G}$ in the Lie algebra of evolution derivations of some 
dif\/ferential ring $C$.  Let $^\sim : C \rightarrow C \otimes {\mathcal G}^*$ 
be the map def\/ined by the relation
\be
\langle H^\sim , X\rangle  \sim X^\wedge (H), \qquad \forall\; H \in C, \quad \forall \;X
\in {\mathcal G}. 
\ee
Extend the given commutator in ${\mathcal G}^*$, -- whether given by an ${\mathcal O}$-operator,
or otherwise, but making ${\mathcal G} + {\mathcal G}^*$ into a dif\/ferential Lie algebra, 
see Appendix~I, -- into the one on $C \otimes 
{\mathcal G}^*$ by treating $C$ as just another dif\/ferential extension $\tilde R \supset R$.  The 
criterion of inf\/initesimal Hamiltonian action then is:
\be
\ba{l}
X^\wedge (\{H, F\}) - \{ X^\wedge (H), F\} - \{ H, X^\wedge (F) \} \sim 
 \langle [H^\sim, F^\sim],  X\rangle, \vspace{2mm}\\
  \forall \; H, F \in C, \quad  \forall\; X \in {\mathcal G}, 
\ea
\ee
with $\{ \ , \}$ being a Poisson bracket on $C$ def\/ined by some Hamiltonian matrix.

We are interested in this paper in the case $C = C_{{\mathcal G}^*} = C_u$.  
In this case the action 
of the evolution vector f\/ield $X^\wedge$ on $C_u$ corresponding to an element 
$X \in {\mathcal G}$  is given by the formula
\be
X^\wedge (\pbf{u}) = X^. u.
\ee

For the {\it linear} bracket on ${\mathcal G}^*$, the RHS of the criterion (5.34) vanishes identically 
since~${\mathcal G}^*$ is considered as an abelian Lie algebra.  Thus, we have to verify that 
\be
X^\wedge (\{H, F\}) \sim \{X^\wedge (H), F\} + \{H, X^\wedge (F)\}, \qquad \forall \; H, F \in C_u. 
\ee
However, this relation follows at once from the fact that $X^\wedge$ is a {\it Hamiltonian} 
vector f\/ield with the Hamiltonian
\be
G = - \langle u, X \rangle.
\ee 
Indeed, by formula (5.5), 
\be
X_G (\pbf{u}) = B \left({\delta G \over \delta \pbf{u}} \right) = - {\delta G \over \delta \pbf{u}} 
^. \pbf{u} = X^. u, 
\ee
and this is formula (5.35) 

We now prove  ``the  main result of the inf\/initesimal Hamiltonian action'': 

\medskip

\noindent
{\bf Theorem 5.39.} {\it The infinitesimal Hamiltonian action criterion (5.34) for a 
given Poisson bracket on ${\mathcal G}^*$ is enough to verify for 
Hamiltonians $H$ and $F$ {\it linear} in the $u$'s. }

\medskip

\noindent
{\bf Proof.} We are going to show that each side of the criterion (5.34) can be 
transformed into a form which is a bilinear {\it differential} operator acting on the vectors
\setcounter{equation}{39}
\be
\pbf{Y} = {\delta H \over \delta \pbf{u}}, \qquad \pbf{Z} = {\delta F \over \delta \pbf{u}}.
\ee

First by formula (5.33), 
\[
\langle H^\sim, X \rangle \sim X^\wedge (H) \sim (X^\wedge (\pbf{u}))^t {\delta H \over \delta \pbf{u}}
 = \langle X^. u, \pbf{Y} \rangle  \sim  \langle u, [Y, X] \rangle  \sim   - \langle Y^. u, X\rangle, 
\]
so that
\be
H^\sim = - {\delta H \over \delta \pbf{u}}^. \pbf{u}.
\ee
Similarly, $F^\sim = - \pbf{Z}^. u,$ and the RHS of the criterion (5.34) is therefore indeed a 
bilinear dif\/ferential operator w.r.t. $\pbf{Y}$ and $\pbf{Z}$.

Next, let $B$ be an unspecif\/ied Hamiltonian matrix over the ring $C_u = C_{{\mathcal G}^{*}}$, so that
\be
\{H, F\} \sim {\delta F \over \delta \pbf{u}^t} B \left({\delta H \over \delta \pbf{u}} \right) 
= \pbf{Z}^t B (\pbf{Y}). 
\ee
Transforming separately each of the 3 terms on the LHS of the criterion (5.34), we get:
\renewcommand{\theequation}{\arabic{section}.\arabic{equation}{\rm a}}
\setcounter{equation}{42}
\be
\ba{ll}
1) & \ds X^\wedge (\{H, F\}) \sim X^\wedge (\pbf{Z}^t B(\pbf{Y}))
\vspace{2mm}\\
& \ds \quad  = (X^\wedge (\pbf{Z}))^t B (\pbf{Y}) + \pbf{Z}^t X^\wedge (B) (\pbf{Y}) +
\pbf{Z}^t B (X^\wedge (\pbf{Y})); 
\ea
\ee
\renewcommand{\theequation}{\arabic{section}.\arabic{equation}{\rm b}}
\setcounter{equation}{42}
\be
\ba{ll}
2) & \ds - \{ X^\wedge (H), F\} \sim - \pbf{Z}^t B {{\delta \over \delta \pbf{u}}}
 (X^\wedge (H))  \vspace{3mm}\\
& \ds \quad  {\mathop{=}\limits^{\mbox{\scriptsize [by  formula (5.46) below]}}} \ 
 - \pbf{Z}^t B ([\pbf{Y}, \pbf{X}] + X^\wedge (\pbf{Y})); 
\ea
\ee
\renewcommand{\theequation}{\arabic{section}.\arabic{equation}{\rm c}}
\setcounter{equation}{42}
\be
\ba{ll}
3) & \ds - \{H, X^\wedge (F)\} \sim \{ X^\wedge (F), H\} \vspace{2mm}\\
& \ds \quad {\mathop{\sim }\limits^{\mbox{\scriptsize  [by (5.43b)]}}} \
 \pbf{Y}^t B ([\pbf{Z}, \pbf{X}] + X^\wedge (\pbf{Z})) \sim - ([ \pbf{Z}, \pbf{X}] + X^\wedge (\pbf{Z}))^t 
B (\pbf{Y}). 
\ea
\ee
Adding up the expressions (5.43), we obtain
\renewcommand{\theequation}{\arabic{section}.\arabic{equation}}
\setcounter{equation}{43} 
\be
\ba{l}
\ds X^\wedge (\{H, F\} ) - \{ X^\wedge (H), F\} - \{H, X^\wedge (F) \} 
\vspace{3mm}\\
\ds \qquad \sim {\delta F \over \delta \pbf{u}^t} X^\wedge (B) 
\left({\delta H \over \delta \pbf{u}} \right) 
+ {\delta F \over \delta \pbf{u}^t} B \left(\left[X, {\delta H \over \delta \pbf{u}}\right]\right) 
+ \left[X, {\delta F \over \delta \pbf{u}} \right]^t B \left({\delta H \over \delta \pbf{u}} \right).
\qquad \mbox{\rule{3mm}{3mm}}
\ea
\ee

\noindent
{\bf Lemma 5.45.}
\setcounter{equation}{45}
\be
{\delta \over \delta \pbf{u}} (X^\wedge (H)) = \left[ {\delta H \over \delta \pbf{u}}, X\right] + 
X^\wedge \left({\delta H \over \delta \pbf{u}} \right). 
\ee

\noindent
{\bf Proof.} We have, 
\be
X^\wedge (H) \sim (X^\wedge (\pbf{u}))^t {\delta H \over \delta \pbf{u}} = 
\langle X ^. u, {\delta H \over \delta \pbf{u}} \rangle  \sim  
\langle u, \left[{\delta H \over \delta \pbf{u}}, X \right] \rangle .
\ee
Therefore, since $X$ is $u$-independent,
\be
{\delta \over \delta \pbf{u}} (X^\wedge (H)) = \left[{\delta H \over \delta \pbf{u}}, X \right] + D 
\left({\delta H \over \delta \pbf{u}}\right)^\dagger (X^. u).
\ee
Now, since
\be
D \left({\delta H \over \delta \pbf{u}} \right)^\dagger = D \left({\delta H \over \delta \pbf{u}} \right)
\ee
(see [6]), the second summand on the RHS of the expression (5.48) can be transformed into 
\[
D \left({\delta H \over \delta \pbf{u}} \right) (X^\wedge (\pbf{u})) = X^\wedge \left({\delta H 
\over \delta \pbf{u}} \right). \hspace{232.2pt} \mbox{\rule{3mm}{3mm}}
\]

\noindent
{\bf Remark 5.50.}  Exactly where have we used the restriction that our objects, 
-- Lie algebras, rings, etc., -- are dif\/ferential rather than dif\/ferential-dif\/ference ones?  
The answer is nowhere except in the notation: in the general case, 
\setcounter{equation}{50}
\be
a \sim b \Leftrightarrow (a - b) \in \sum^m_{\ell =1} \mbox{Im}\, \partial_\ell + \sum_{g \in G} 
\mbox{Im}\, (\hat g - \hat e),
\ee
$C_u = R [u_i^{(g|\sigma)}]$, $g \in G$, $\sigma \in {\mathbf Z}^m_+$, etc., 
where $G$ is a discrete 
group whose elements index the discrete degrees of freedom, $\hat g$ is the action of the 
element $g \in G$ on $R$, $C_u$, etc.  The presence of discrete degrees of freedom is hidden 
in the notation ${{\ds{\delta \over \delta \pbf{u}}}}$, $( \ )^\dagger$, $D ( \cdot )$, etc.  
(The continuous reader may safely ignore this Remark).

\medskip

We conclude this Section by considering {\it affine} Poisson brackets.  These brackets have the 
corresponding Hamiltonian operators of the form
\be
B = B^{\mbox{\scriptsize lin}} + b,
\ee
where $B^{\mbox{\scriptsize lin}}$ is a Hamiltonian matrix operator {\it linearly} 
dependent upon $\pbf{u}$, and 
$b$ is $\pbf{u}$-independent.  Thus,
\be
\{H, F\} \sim {\delta F \over \delta \pbf{u}^t} B \left({\delta H \over \delta \pbf{u}}\right) 
= \langle \pbf{u}, \left[{\delta H \over \delta \pbf{u}}, {\delta F \over \delta \pbf{u}}\right] \rangle
 + \langle b \left({\delta
 H \over \delta \pbf{u}} \right), {\delta F \over \delta \pbf{u}} \rangle, 
\ee
where $[\ , ]$ is a commutator in some Lie algebra, say ${\mathcal G}$, and 
$b: {\mathcal G} \rightarrow {\mathcal G}^*$ is 
a skewsymmetric operator def\/ining a generalized 2-cocycle on ${\mathcal G}$:
\be
\langle b ([X, Y]), Z \rangle + \mbox{c.p.} \sim 0, \qquad  \forall\; X, Y, Z \in {\mathcal G}
\ee
(see [6].)  Since thie 2-cocycle is, in general, {\it generalized} (``$\sim$'' instead of 
``='' in the RHS of (5.54)), it does not correspond to a central extension any more.  
Nevertheless, we have

\medskip

\noindent
{\bf Theorem 5.55.} {\it  Let ${\mathcal G}$ act on $C_{{\mathcal G}^{*}}$ by the rule
\setcounter{equation}{55}
\be
X^\wedge (\pbf{u}) = X ^. u - b (\pbf{X}). 
\ee
Then this action satisfies the infinitesimal Hamiltonian action criterion (5.34) for the 
affine Poisson bracket (5.53).}

\medskip

\noindent
{\bf Proof.} Let us write
\be
X^\wedge = X^\wedge_{\mbox{\scriptsize old}} + X^\wedge_{\mbox{\scriptsize new}}, 
\qquad  \{H, F\} = \{H, F\}_{\mbox{\scriptsize old}} + \{H, F\}_{\mbox{\scriptsize new}}, 
\ee
where
\be
X^\wedge_{\mbox{\scriptsize new}} (\pbf{u}) = - b (\pbf{X}),
\ee
\be
\{H, F\}_{\mbox{\scriptsize new}} = \langle b \left({\delta H \over \delta \pbf{u}} \right), 
{\delta F \over \delta \pbf{u}} \rangle . 
\ee
By Theorem 5.39, we need to verify the relation
\be
X^\wedge (\{H, F\}) - \{X^\wedge (H), F\} - \{H, X^\wedge (F) \} \sim 0
\ee
for all $H$, $F$ linear in $\pbf{u}$:
\be
H = \langle u, Y\rangle, \qquad  F = \langle u, Z \rangle, \qquad \forall\; Y, Z \in {\mathcal G}.
\ee
We have:
\renewcommand{\theequation}{\arabic{section}.\arabic{equation}{\rm a}}
\setcounter{equation}{61}
\be
\ba{l}
\ds X^\wedge (\{H, F\}) \sim (X^\wedge_{\mbox{\scriptsize old}} + 
X^\wedge_{\mbox{\scriptsize new}} ) \{ \langle u, [Y, Z] \rangle + \langle b (Y), Z \rangle \}
\vspace{2mm}\\
\ds \qquad = X^\wedge_{\mbox{\scriptsize old}} (\{H, F\}_{\mbox{\scriptsize old}}) 
- \langle b (X), [Y, Z] \rangle,
\ea
\ee
\renewcommand{\theequation}{\arabic{section}.\arabic{equation}{\rm b}}
\setcounter{equation}{61}
\be
\ba{l}
\ds - \{X^\wedge(H), F\} = - \{ \langle X^. \ u - b (X), Y \rangle, \langle u, Z \rangle \} 
\vspace{2mm}\\
\ds \qquad \sim \ \{ \langle u, [X, Y] \rangle , \langle u, Z \rangle \} = 
- \{X^{\mbox{\scriptsize old}} (H), F\}_{\mbox{\scriptsize old}} + \langle b( [X, Y]), Z \rangle , 
\ea
\ee
\renewcommand{\theequation}{\arabic{section}.\arabic{equation}{\rm c}}
\setcounter{equation}{61}
\be
- \{H, X^\wedge (F) \} \sim \ - \{ H, X^\wedge_{\mbox{\scriptsize old}} (F)\}_{\mbox{\scriptsize old}}
 - \langle b ([X, Z]), Y\rangle . 
\ee
Adding the expressions (5.62) up and remembering formula (5.36), we obtain 
formula (5.60).~\hfill\rule{3mm}{3mm}

\renewcommand{\theequation}{\arabic{section}.\arabic{equation}}
\setcounter{equation}{0}

\section{ Quadratic Poisson brackets on dual spaces to Lie algebras}

The action map $\mbox{Ad}^*: G \times {\mathcal G}^* \rightarrow {\mathcal G}^*$ is 
Poisson when the Poisson bracket on 
$G$ is zero and on ${\mathcal G}^*$ is linear.  When $G$ itself has 
a nonzero multiplicative Poisson 
bracket on it, {\it coming from} an $r$-matrix, there exists a quadratic deformation of the 
standard linear Poisson bracket on ${\mathcal G}^*$ such that the action map 
$\mbox{Ad}^*: G \times {\mathcal G}^* \rightarrow 
{\mathcal G}^*$ is still Poisson.  (If the multiplicative Poisson bracket on $G$ is {\it not} of 
$r$-matrix type, such quadratic deformation is, in general, impossible.)  This was found by 
Kupershmidt and Stoyanov in~[11].  In this Section we construct a dif\/ferential (-dif\/ference) 
analog of this quadratic bracket, prove that it's compatible with the linear one, verify that 
this quadratic bracket is natural, and then check the inf\/initesimal Hamiltonian action 
criterion for it.  

So, let ${\mathcal O}: {\mathcal G}^* \rightarrow {\mathcal G}$ be an ${\mathcal O}$-operator.  
Recall that ${\mathcal O}$ is skewsymmetric: 
\be
\langle u, {\mathcal O} (v) \rangle  \sim - \langle v, {\mathcal O} (u) \rangle, \qquad
 \forall \; u, v \in {\mathcal G}^*, 
\ee
and
\be
{\mathcal O} ({\mathcal O} (u)^. v - {\mathcal O} (v)^. u) = [{\mathcal O} (u), {\mathcal O}(v)], 
\qquad  \forall \; u, v \in {\mathcal G}^*. 
\ee

Def\/ine the quadratic Poisson bracket on ${\mathcal G}^*$ by formula
\renewcommand{\theequation}{\arabic{section}.\arabic{equation}{\rm a}}
\setcounter{equation}{2}
\be
\{H, F\} \sim \langle {\delta H \over \delta \pbf{u}}^. \pbf{u}, {\mathcal O} \left({\delta F 
\over \delta \pbf{u}}^. u \right) \rangle  \sim - \langle {\delta F  \over \delta \pbf{u}}^.  \pbf{u}, 
{\mathcal O} \left({\delta H \over \delta \pbf{u}}^. \pbf{u} \right) \rangle
\ee
\renewcommand{\theequation}{\arabic{section}.\arabic{equation}{\rm b}}
\setcounter{equation}{2}
\be
\qquad \sim \langle {\mathcal O} \left({\delta H \  \over \delta \pbf{u}}^. \pbf{u} \right)^. 
\pbf{u}, {\delta F \over \delta \pbf{u}} \rangle . 
\ee
The corresponding Hamiltonian matrix $B$ is therefore quadratic in $u$:
\renewcommand{\theequation}{\arabic{section}.\arabic{equation}}
\setcounter{equation}{3}
\be
B (\pbf{Y}) = {\mathcal O} (\pbf{Y}^. \pbf{u})^. \pbf{u}. 
\ee

The quadratic bracket (6.3) is obviously skewsymmetric.  Let us verify that it satisf\/ies the 
Jacobi identity.  By the main result of the Hamiltonian formalism ([6]~p.~47), it is enough 
to check the Jacobi identify for Hamiltonians {\it linear} in $\pbf{u}$.  So, let 
\be
H = \pbf{Y}^t \pbf{u}, \qquad  F = \pbf{Z}^t \pbf{u}, \qquad G = \pbf{X}^t \pbf{u}. 
\ee
Then
\be
\{H, F\} \sim \langle \pbf{Y}^. \pbf{u}, {\mathcal O} (Z^. \pbf{u}) \rangle 
\sim - \langle \pbf{u}, [\pbf{Y}, {\mathcal O} (\pbf{Z}^. \pbf{u})] \rangle
\sim  \langle \pbf{u}, [\pbf{Z}, {\mathcal O} (\pbf{Y}^. \pbf{u})] \rangle .
\ee
Therefore,
\be
{\delta \{H, F \} \over \delta \pbf{u}} = - [\pbf{Y}, {\mathcal O} (\pbf{Z}^. \pbf{u}) ] + 
[ \pbf{Z}, {\mathcal O} (\pbf{Y}^.\pbf{u}) ], 
\ee
and hence, in the notation
\be
\bar X = X^. u, \qquad  \bar Y = Y^. u, \qquad  \bar Z = Z^. u,
\ee
\[
\ba{l}
\ds \{\{ H, F \}, G \} + \mbox{c.p.} \sim \langle - [ Y, {\mathcal O} (\bar Z)]^.u 
+ [Z, {\mathcal O} (\bar Y) ]^. u, {\mathcal O} (\bar X) \rangle + \mbox{c.p.}
\vspace{2mm}\\
\ds \qquad = - (\langle[Y, {\mathcal O} (\bar Z)]^. u, {\mathcal O} (\bar X) \rangle + \mbox{c.p.}) 
+ (\langle [Y, {\mathcal O} (\bar X)]^. u, {\mathcal O} (\bar Z) \rangle + \mbox{c.p.})
\vspace{2mm}\\
\ds \qquad \sim  \langle u, [[Y, {\mathcal O} (\bar Z)], {\mathcal O} (\bar X)] 
- [[Y, {\mathcal O} (\bar X)], {\mathcal O} (\bar Z)] \rangle + \mbox{c.p.}
\vspace{2mm}\\
\ds \qquad = \langle u, - [Y, [{\mathcal O} (\bar X), {\mathcal O} (\bar Z)] \rangle + 
\mbox{c.p.} \sim  \langle \bar Y, [ {\mathcal O} (\bar X), {\mathcal O} 
(\bar Z)]\rangle + \mbox{c.p.}
\ea
\]
and this expression is $\sim$ 0 by formula (2.6), itself equivalent to the
${\mathcal O}$-property~(6.2).

Now let us verify that the linear and quadratic Poisson brackets on ${\mathcal G}^*$ 
are compatible no 
matter what ${\mathcal O}$ is.  Recall that compatibility of two Poisson brackets means that their 
arbitrary linear combination with constant coef\/f\/icients is again a Poisson bracket, i.e., it 
satisf\/ies the Jacobi identify.  This amounts to the relation
\be
\left( \{ \{ H, F \}_1, G\}_2 + \{ \{ H, F\}_2, G \}_1 \right)+ \mbox{c.p.} \sim 0, \qquad 
\forall\; H, G, F, 
\ee
and the main Theorem of the Hamiltonian formalism asserts that this relation needs to be 
verif\/ied only for {\it linear} Hamiltonian $H$, $F$, $G$.  So, for such $H$, $F$, $G$, given by 
formula (6.5), we have, by formulae (5.2) and (6.6) 
\renewcommand{\theequation}{\arabic{section}.\arabic{equation}{\rm a}}
\setcounter{equation}{9}
\be
\{ H, F\}_1 \sim [\pbf{Y}, \pbf{Z}]^t \pbf{u} \ \Rightarrow \
\{ \{ H, F \}_1, G \}_2 \sim \langle [ \pbf{Y}, \pbf{Z}]^. \pbf{u}, {\mathcal O} (\pbf{X}^. \pbf{u}) \rangle.
\ee
On the other hand, by formulae (6.7) and (5.6), we get 
\renewcommand{\theequation}{\arabic{section}.\arabic{equation}{\rm b}}
\setcounter{equation}{9}
\be
\ba{l}
\ds \{ \{H, F\}_2, G\}_1 \sim \langle \pbf{X}^. \pbf{u}, - [\pbf{Y}, {\mathcal O}(\pbf{Z}^. u) ] 
+ [\pbf{Z}, {\mathcal O} (\pbf{Y}^. u)] \rangle
\vspace{2mm}\\
\ds \qquad \sim \langle  \pbf{Y}^. (\pbf{X}^. \pbf{u}), {\mathcal O} (\pbf{Z}^. u) \rangle
 - \langle \pbf{Z}^. (\pbf{X}^. \pbf{u}), {\mathcal O} (\pbf{Y}^. \pbf{u}) \rangle.
\ea
\ee
Substituting expressions (6.10) into formula (6.9), we get
\[
(\langle [\pbf{X}, \pbf{Y}]^. \pbf{u} + \pbf{Y} (\pbf{X}^. \pbf{u})
 - \pbf{X}^. (\pbf{Y}^. \pbf{u}), {\mathcal O} (\pbf{Z}^. u) \rangle  ) + \mbox{c.p.} = 0 ,
\]
since
\renewcommand{\theequation}{\arabic{section}.\arabic{equation}}
\setcounter{equation}{10}
\be
[\pbf{X}, \pbf{Y}]^. \pbf{u} = \pbf{X}^. (\pbf{Y}^. \pbf{u}) - \pbf{Y}^. (\pbf{X}^. \pbf{u}). 
\ee

If it so happens that the ${\mathcal O}$-operator ${\mathcal O}: 
{\mathcal G}^* \rightarrow {\mathcal G}$ is invertible, like in \S~2, 
then we have a generalized 2-cocycle on ${\mathcal G}$:
\renewcommand{\theequation}{\arabic{section}.\arabic{equation}{\rm a}}
\setcounter{equation}{10}
\be
\omega_b (X, Y) = \langle b (X),  Y\rangle, \qquad
 b = \epsilon {\mathcal O}^{-1},  \qquad  \epsilon = \mbox{const}, 
\ee
and thus a constant-coef\/f\/icient Poisson bracket on $C_{{\mathcal G}{^*}}$: 
\renewcommand{\theequation}{\arabic{section}.\arabic{equation}{\rm b}}
\setcounter{equation}{10}
\be
\{ H, F\}_0 \sim \langle b \left({\delta H \over \delta u} \right), {\delta F \over \delta u} \rangle. 
\ee
Since $\omega_b$ is a generalized 2-cocycle on ${\mathcal G}$, 
this constant-coef\/f\/icient Poisson bracket 
$\{\ , \}_0$ on ${\mathcal G}^*$ is compatible with the linear Poisson bracket $\{\ , \}_1$.  
Let us verify 
that all three Poisson brackets on ${\mathcal G}^*$, -- constant-coef\/f\/icient, linear and 
quadratic, -- are 
compatible.  It remains only to verify compatibility of constant-coef\/f\/icient $\{\ , \}_0$ one 
and the quadratic  $\{\ , \}_2$ one.  Again, for linear Hamiltonians (6.5), the Poisson bracket 
$\{H, F\}_0$ (6.11b) is $u$-independent, so that $\{\{H, F\}_0, G\}_2 = 0$.  Thus, we need 
only to verify that 
\[
\{\{H, F\}_2, G\}_0 + c.p. \sim 0.
\]
By formulae (6.7), (6.8), and (6.11), we have:
\[
\ba{l}
\ds \{\{H, F\}_2, G\}_0 + \mbox{c.p.} = \langle b \left({\delta \{ H, F\}_2 \over \delta u}\right), X 
\rangle + \mbox{c.p.} 
\vspace{3mm}\\
\ds \qquad \sim \langle b (X), [Y, {\mathcal O} (\bar Z] - [Z, {\mathcal O} (\bar Y)] \rangle 
+ \mbox{c.p.}
\vspace{2mm}\\
\ds \qquad = ( \langle b (X), [Y, {\mathcal O} (\bar Z)] \rangle - \langle b (Y), [X, {\mathcal O} (Z)] 
\rangle ) + \mbox{c.p.} 
\vspace{2mm}\\
\ds \qquad \sim \langle - Y^. b (X) + X^. b (Y), {\mathcal O} (\bar Z) \rangle + \mbox{c.p.}
\vspace{2mm}\\
\ds \qquad {\mathop{=}\limits^{\mbox{\scriptsize [by (2.11)]}}} \
 \langle \epsilon [{\mathcal O}^{-1} (X), {\mathcal O}^{-1} (Y)],  {\mathcal O} (\bar Z) \rangle + 
\mbox{c.p.} 
\vspace{2mm}\\
\ds \qquad \sim - \epsilon \langle Z^. u,  {\mathcal O} [{\mathcal O}^{-1} (X), 
{\mathcal O}^{-1} (Y) ] \rangle + \mbox{c.p.}
\vspace{2mm}\\
\ds \qquad {\mathop{=}\limits^{\mbox{\scriptsize [by (2.12)]}}} \
 - \epsilon \langle  Z^. u,  [ X, Y] \rangle + \mbox{c.p.}  \sim \epsilon 
\langle u, [Z, [X, Y]] \rangle + \mbox{c.p.} = 0 .
\ea
\]

Thus, when ${\mathcal O}$ is invertible, we have a triple of compatible Hamiltonian structures on 
${\mathcal G}^*$, of $u$-degrees zero, one, and two.  When ${\mathcal O}$ is not invertible, we are left with 
only linear and quadratic Poisson brackets.

To see that the quadratic Poisson bracket is {\it{natural}}, let $\varphi: {\mathcal G} \rightarrow 
{\mathcal H}$ be a homomorphism of Lie algebras.  
Let ${\mathcal O}_{{\mathcal H}} = \varphi {\mathcal O}_{{\mathcal G}} \varphi^\dagger: 
{\mathcal H}^* \rightarrow {\mathcal H}$ be the ${\mathcal O}$-operator induced 
on ${\mathcal H}^*$ by the ${\mathcal O}$-operator 
${\mathcal O} = {\mathcal O}_{{\mathcal G}}$ on ${\mathcal G}^*$.  
Let $\Phi:  C_u \rightarrow C_q$ (5.20) be the corresponding 
homomorphism of function rings:
\renewcommand{\theequation}{\arabic{section}.\arabic{equation}}
\setcounter{equation}{11}
\be
\Phi (\pbf{u}) = \varphi^\dagger (\pbf{q}). 
\ee
To show that the map $\Phi$ is Hamiltonian between the quadratic Poisson brackets on 
${\mathcal G}^*$ and ${\mathcal H}^*$, we use Proposition~5.24.  
So, let $H$ and $F$ be two linear Hamiltonians in $C_u$:
\be
H = \pbf{Y}^t \pbf{u}, \qquad F = \pbf{Z}^t \pbf{u}, \qquad 
\pbf{Y}, \pbf{Z} \in {\mathcal G}. 
\ee
Then, by formulae (6.3) and (6.12), 
\be
\ba{l}
\Phi (\{H, F\}_{{\mathcal G}^{*}}) \sim \Phi (\langle
 \pbf{Y}^. \pbf{u}, {\mathcal O}_{\mathcal G} (\pbf{Z}^. \pbf{u}) \rangle ) 
= \langle \pbf{Y}^. (\varphi^\dagger (\pbf{q}), {\mathcal O}_{{\mathcal G}} 
(Z^. \varphi^\dagger (\pbf{q})) \rangle
\vspace{2mm}\\
\ds \qquad {\mathop{=}\limits^{\mbox{\scriptsize [by (4.15)]}}} \
 \langle \varphi^\dagger (\varphi (\pbf{Y})^. \pbf{q}), {\mathcal O}_{\mathcal G} 
\varphi^\dagger (\varphi (Z)^. \pbf{q})  \rangle  \sim  
\langle \varphi (\pbf{Y})^. \pbf{q}, \varphi {\mathcal O}_{\mathcal G} 
\varphi^\dagger (\varphi (\pbf{Z})^. \pbf{q}) \rangle
\vspace{3mm}\\
\ds \qquad = \langle \varphi (\pbf{Y})^. \pbf{q}, {\mathcal O}_{{\mathcal H}} (\varphi (\pbf{Z})^. \pbf{q}) 
\rangle = \langle {\delta \Phi (H) \over 
\delta \pbf{q}}^. \pbf{q}, {\mathcal O}_{{\mathcal H}} \left( {\delta \Phi (F) \over \delta \pbf{q}}^. 
\pbf{q} \right) \rangle \sim  \{ \Phi (H), \Phi (F) \}_{{\mathcal H}^{*}}, 
\ea\hspace{-19.07pt}
\ee
where we used in the next to last equality in this chain the formula
\[
{\delta \Phi (H) \over \delta \pbf{q}} = \varphi (\pbf{Y}), 
\]
which follows from the relations
\[
\Phi (H) = \Phi (\langle \pbf{u}, \pbf{Y})\rangle ) = \langle \varphi^\dagger (\pbf{q}), \pbf{Y}) 
\rangle \sim \langle \pbf{q}, \varphi (\pbf{Y}) \rangle. 
\]

What are the Casimirs of the quadratic bracket?  Formula (6.4) shows that they are precisely the  
solutions of the equation
\be
{\mathcal O} \left( {\delta H \over \delta \pbf{u}}^. \pbf{u} \right)^. \pbf{u} = 0. 
\ee
In particular, all ``coadjoint invariants'', i.e., those $H$ satisfying
\be
{\delta H \over \delta \pbf{u}}^. \pbf{u} = 0, 
\ee
are Casimirs, so that, in f\/inite dimensions, the symplectic leaves of the quadratic bracket sit 
inside the coadjoint orbits.  A better understanding of the symplectic leaves should be interesting.  

Let us now check the inf\/initesimal Hamiltonian action criterion for the quadratic bracket.  
By Theorem 5.39, we have to verify the relation (5.34) for linear Hamiltonians $H$ and $F$ 
given by formula (6.13).  Starting with the RHS of the criterion (5.34) and using formulae 
(5.41) and (2.11), we obtain:
\be
\ba{l}
\ds \langle [H^\sim, F^\sim], X \rangle = \langle [- Y^. u, - Z^. u], X\rangle 
\vspace{2mm}\\
\ds \qquad = \langle {\mathcal O} (Y^. u)^. (Z^. u) - {\mathcal O} (Z^. u)^. (Y^. u), X \rangle
 = \langle {\mathcal O} (\bar Y)^. \bar Z - {\mathcal O} (\bar Z)^. \bar Y, X \rangle, 
\ea
\ee
where we introduced the convenient notation
\be
\bar Y = Y^. u, \qquad  \bar Z = Z^. u. 
\ee
For LHS of the criterion (5.34) we use the form (5.44):
\be
Z^t (X ^\wedge (B) (Y) + Z^t B ([X, Y]) + [X, Z]^t B(Y)
\ee
\[
 {\mathop{\sim}\limits^{\mbox{\scriptsize [by (6.4)]}}} \ 
Z^t ({\mathcal O} (Y^. (X^. u))^. u + {\mathcal O} (Y^. u) ^. (X^. u)) - [X, Y]^t {\mathcal O} (Z^. u)^. 
u +  [X, Z]^t {\mathcal O} (Y^. u)^. u 
\]
\be
 {\mathop{\sim}\limits^{\mbox{\scriptsize [by (6.18), (6.22)]}}} \ 
 \langle - Z^. u, {\mathcal O} (Y^. (X^. u)) \rangle + \langle X^. u, [Z, {\mathcal O} (\bar Y)]
\rangle
\ee
\be
 + \langle u, [ {\mathcal O} (\bar Z), [ X, Y]] \rangle - \langle u, [ {\mathcal O} 
(\bar Y), [ X, Z]] \rangle .
\ee
(We used above the obvious relation
\be
X^t (Y^. u) \sim - Y^t (X^. u) . ) 
\ee
The 1$^{st}$ summand in the expression (6.20) can be transformed as
\renewcommand{\theequation}{\arabic{section}.\arabic{equation}{\rm a}}
\setcounter{equation}{22}
\be
\langle
 Y^. (X^. u), {\mathcal O} (\bar Z) \rangle \sim \langle u, [X, [Y, {\mathcal O} (\bar Z)]]\rangle,
\ee
while the second summand in (6.20) is $\sim$ to 
\renewcommand{\theequation}{\arabic{section}.\arabic{equation}{\rm b}}
\setcounter{equation}{22}
\be
- \langle u, [X, [Z, {\mathcal O} (\bar Y)]]\rangle.
\ee
Altogether, expression (6.21) and (6.23) add up to
\renewcommand{\theequation}{\arabic{section}.\arabic{equation}}
\setcounter{equation}{23}
\be
\langle
 u, - [Y, [ {\mathcal O} (\bar Z), X]] + [Z, [ {\mathcal O} (\bar Y), X]] \rangle
\sim \langle {\mathcal O} (\bar Z)^. (Y^. u) - {\mathcal O} (\bar Y)^. (Z^. u), X \rangle,
\ee
and this is the same as the expression (6.17).

\medskip

\noindent
{\bf Example 6.25.}  Let ${\mathcal G} = {\mathcal G} (\mu)$ be the Lie algebra (3.7), 
\setcounter{equation}{25}
\be
\left[{X \choose f}, {Y \choose g} \right] =  {XY^\prime - X^\prime Y \choose (Xg - Yf + 
\mu (X^\prime Y^{\prime \prime} - X^{\prime \prime} Y^\prime) )^\prime} ,
\ee
and let ${\mathcal O}$ be the ${\mathcal O}$-operator (3.12):
\be
{\mathcal O} = \left(\matrix{0 & 1 \cr -1 & \epsilon \partial^3 \cr} \right). 
\ee
By formula (3.14),
\be
{X \choose f}^. {u \choose p} = {Xu^\prime + 2X^\prime u - f p^\prime + \mu (X^\prime p^\prime 
)^{\prime \prime} + \mu (X^{\prime \prime} p^\prime)^\prime \choose Xp^\prime} . 
\ee
For $H \in C_{u, p},$ let
\be
{X \choose f} = {\delta H/\delta u \choose \delta H/\delta p}. 
\ee
Then the motion equations for the Hamiltonian vector f\/ield $X_H$, by formula (6.4), are
\be
\ba{l}
\ds {u \choose p}_t = X_H {u \choose p} = B (\pbf{X}) = 
B {\delta H/\delta u \choose \delta H/\delta 
p} = {\mathcal O} \bigg(\pbf{X}^. {u \choose p} \bigg)^. {u \choose p} 
\vspace{4mm}\\
\ds \qquad = \left(\left(\matrix{0 & 1 \cr -1 & \epsilon \partial^3 \cr} \right)  \left(\matrix{
Xu^\prime + 2 X^\prime u - f p^\prime + \mu (X^\prime p^\prime ) ^{\prime \prime} + \mu 
(X^{\prime \prime} p^\prime)^\prime \cr Xp^\prime \cr} \right) \right)^. {u \choose p} 
\vspace{4mm}\\
\ds \qquad = \left(\matrix{Xp^\prime \cr \epsilon (Xp^\prime)^{\prime \prime \prime} - Xu^\prime - 
2X^\prime u - \mu (X^\prime p^\prime)^{\prime \prime} - \mu (X^{\prime \prime} p^\prime)^\prime 
+ f p^\prime \cr} \right)^. {u \choose p} 
\vspace{4mm}\\
\ds \qquad = \left(\matrix{Xp^\prime u^\prime + 2 (Xp^\prime)^\prime u + p^\prime (- \epsilon (Xp^\prime) 
^{\prime \prime \prime} - f p^\prime + Xu^\prime + 2 X^\prime u + \mu (X^\prime p^\prime)^
{\prime \prime}  \cr 
+ \mu (X^{\prime \prime} p^\prime)^\prime + \mu ( (Xp^\prime)^\prime p^ \prime)^{\prime \prime}+ \mu ((Xp^\prime)^{\prime \prime} p^\prime)
^\prime \cr
- - - - - - - - - - - - - - - - - - - - - - - - - - - -  \cr
Xp^{\prime 2} \cr} \right).
\ea\!\!
\ee
Hence, the quadratic Hamiltonian matrix on ${\mathcal G}(\mu^*$ is
\be
B = \left(\matrix{ * & -p^{\prime  2} \cr p^{\prime 2} & 0 \cr} \right), 
\ee
\be
* = 2 (p^\prime u \partial + \partial p^\prime u) + (4 \mu - \epsilon) p^\prime \partial^3 
p^\prime + \mu (3 p^{\prime \prime 2} - 2 p^\prime p^{\prime \prime \prime} ) 
\partial + \mu \partial (3 p^{\prime \prime 2} - 2 p^\prime p^{\prime \prime \prime}). 
\ee
This quadratic Hamiltonian matrix is compatible with the linear Hamiltonian matrix
\renewcommand{\theequation}{\arabic{section}.\arabic{equation}{\rm a}}
\setcounter{equation}{32}
\be
B^{\mbox{\scriptsize lin}} = - \left(\matrix{** & -p^\prime \cr p^\prime & 0 \cr}\right), 
\ee
\renewcommand{\theequation}{\arabic{section}.\arabic{equation}{\rm b}}
\setcounter{equation}{32}
\be
** = u \partial + \partial u + \mu \left(\partial^2 p^\prime \partial + \partial p^\prime 
\partial^2\right), 
\ee
and both these Hamiltonian matrices are compatible with the constant-coef\/f\/icient Hamiltonian 
matrix
\renewcommand{\theequation}{\arabic{section}.\arabic{equation}}
\setcounter{equation}{33}
\be
{\mathcal O}^{-1} = \left(\matrix{\epsilon \partial^3 & -1 \cr 1 & 0 \cr} \right). 
\ee
For $\mu = \epsilon = 0$, formulate (6.31)--(6.34) have $n$-dimensional analogs. 
See Appendix~A2. 

We conclude this Section by calculating the quadratic Poisson bracket in f\/inite dimensions.  
Let $(e_i) $ be a basis in ${\mathcal G}, (e^i)$ the dual basis in ${\mathcal G}^*$,  
$r = \sum\limits_{ij} r^{ij} 
e_i \otimes e_j \in {\mathcal G}^{\otimes 2}$, $r^{ij} = - r^{ji}$, the classical $r$-matrix, 
$(c_{ij}^k )$ the structure constants of ${\mathcal G}$ in the chosen basis.  Then
\be
{\mathcal O} (e^i ) = \sum_s e_s r^{si}, 
\ee
\be
e_j ^\cdot  e^j = - \sum_s c^j_{is} e^s, 
\ee
\be
\ba{l}
\ds \{u_i, u_j\} = \langle e_i^\cdot \pbf{u}, {\mathcal O} (e_j^\cdot \pbf{u}) \rangle 
= \sum_{\alpha \beta} u_\alpha u_\beta 
\langle e_i ^\cdot  e^\alpha, {\mathcal O} (e_j^\cdot  e_\beta) \rangle
\vspace{3mm}\\
\ds \qquad = \sum_{\alpha \beta k \ell} u_\alpha u_\beta c^\alpha_{ik} c^\beta_{j\ell} 
\langle  e^k, {\mathcal O} (e^\ell) \rangle = \sum_{\alpha \beta k \ell} u_\alpha u_\beta 
c^\alpha_{ik} c^\beta_{j\ell} r^{k \ell}. 
\ea
\ee
This is formula (28) in [11].  All other results of this Section had been established for the 
f\/inite-dimensional case in that paper.

\setcounter{equation}{0}

\section{Symplectic models for linear Poisson brackets\\
\hspace*{12mm}on dual spaces to Lie algebras}

Let $\chi: {\mathcal G} \rightarrow \mbox{Dif\/f}\, (V)$ be a representation of a Lie 
algebra ${\mathcal G}$ on a vector space $V$.  Let
\be
\nabla: \ V \times V^* \rightarrow {\mathcal G}^* 
\ee
be the map def\/ined by the relation
\be
\langle v \nabla v^*, X \rangle \sim \langle v^*, \chi (X) (v) \rangle, 
\qquad  \forall \; v \in V, \ v^* \in V^*, \ X \in {\mathcal G}. 
\ee
This map is then {\it Hamiltonian}, between the linear Poisson bracket on ${\mathcal G}^*$ and 
symplectic Poisson bracket on $V \oplus V^*$.  This was proven in~[5] Ch.~8, where I called 
such maps Clebsch representations.  We shall see in the next Section 
that the {\it same  map} (7.1) is Hamiltonian between the quadratic Poisson bracket 
on ${\mathcal G}^*$ 
def\/ined in \S~6 and some interesting {\it quadratic} Poisson bracket on $V \oplus V^*$.  
In this Section we prepare the ground for the next one, by f\/ixing notation and quickly 
reproving the Hamiltonian property of the map $\nabla$ (7.1) for the linear Poisson bracket.

Let $C_{\mathcal M} = \mbox{Fun}\, (V \oplus V^*) = R[ x_\alpha^{(\sigma)}, 
p_\alpha^{(\sigma)}]$, $\alpha = 1, \ldots, \dim (V)$.  
Def\/ine the symplectic Poisson bracket on $C_{\mathcal M}$ by the matrix
\be
b = \left(\matrix{{\mathbf 0} & - {\mathbf 1} \cr {\mathbf 1} & {\mathbf 0} \cr} \right),
\ee
so that
\be
\{H, F\} \sim \left(\begin{array}{c} {\delta F \over \delta \pbf{x}} \vspace{2mm}\\ 
{\delta F \over 
\delta \pbf{p}} \end{array}
 \right)^t  b \left(\begin{array}{c}{\delta H \over \delta \pbf{x}} \vspace{2mm}\\
{\delta H \over \delta \pbf{p}} \end{array} \right) = - {\delta F \over \delta \pbf{x}^t} {\delta H \over 
\delta \pbf{p}} + { \delta F \over \delta \pbf{p}^t} {\delta H \over \delta \pbf{x}}. 
\ee
Let $\Phi:  C_{\mathcal G}^* \rightarrow C_{\mathcal M}$ be the dif\/ferential (-dif\/ference) 
homomorphism def\/ined on the generators of the ring $C_{{\mathcal G}^*} = C_u = 
R[u_i^{(\sigma)}]$, $i = 1, \ldots , \dim ({\mathcal G})$, by the rule
\be
\Phi (\pbf{u}) = \pbf{x} \nabla \pbf{p}.
\ee
To show that this map $\Phi$ is Hamiltonian, we appeal to Propositions 5.24, choose two linear 
in $u$ Hamiltonians
\be
H = \pbf{u}^t \pbf{Y}, \qquad F = \pbf{u}^t \pbf{Z}, \qquad \pbf{Y}, \pbf{Z} \in {\mathcal G}, 
\ee
and then have:
\be
\Phi (\{H, F\}_{{\mathcal G}^{*}} ) \sim \Phi (\pbf{u}^t [\pbf{Y}, \pbf{Z}]) = \langle
\pbf{x} \nabla \pbf{p}, [ \pbf{Y}, \pbf{Z}] \rangle  \sim
 \langle \pbf{p}, [ \pbf{Y}, \pbf{Z}]. \pbf{x} \rangle ,
\ee
\be
\Phi (H) = \Phi (\pbf{u}^t \pbf{Y}) = \langle \pbf{x} \nabla \pbf{p}, \pbf{Y} \rangle \sim 
\langle \pbf{p}, \pbf{Y}. \pbf{x} \rangle \sim \langle - \pbf{Y}^. \pbf{p}, \pbf{x} \rangle
\ee
\be  
\Rightarrow  \ {\delta \Phi (H) \over \delta \pbf{x}} = - \pbf{Y}^. \pbf{p}, \qquad
 {\delta \Phi (H) \over \delta \pbf{p}} = \pbf{Y}. \pbf{x} 
\ee
\be
\ba{l}
\Rightarrow \ 
\{\Phi (H), \Phi (F) \}_{\mathcal M} \ {\mathop{\sim}\limits^{\mbox{\scriptsize [by (7.4)]}}}
\  (\pbf{Z}^. \pbf{p})^t (\pbf{Y} . \pbf{x}) - (\pbf{Z} . \pbf{x})^t (\pbf{Y}^. \pbf{p}) 
\vspace{2mm}\\
\qquad \sim \langle \pbf{p}, - \pbf{Z}. (\pbf{Y}. \pbf{x}) \rangle + 
\langle \pbf{p}, \pbf{Y}. (\pbf{Z}. \pbf{x}) \rangle = \langle \pbf{p}, - \pbf{Z}. 
(\pbf{Y}. \pbf{x}) + \pbf{Y}. (\pbf{Z}. \pbf{x}) \rangle
\vspace{2mm}\\
\qquad = \langle \pbf{p}, [\pbf{Y}, \pbf{Z}]. \pbf{x} \rangle, 
\ea
\ee
and this is the same as the expression (7.7).

\medskip

\noindent
{\bf Remark 7.11.} Strictly speaking, we are dealing here not with the {\it full} 
Clebsch representations, -- which are Hamiltonian maps:  
$\mbox{Fun}\, (({\mathcal G} {\mathop{\ltimes}\limits_{\chi}} V^*)
\rightarrow \mbox{Fun}\, (V \oplus V^*)$, -- but only with their ${\mathcal G}^*$-components, in which the nontrivality 
of results resides; and in any case, this component is the only one we need in this paper.

\setcounter{equation}{0}

\section{Clebsch representations for quadratic Poisson brackets\\
\hspace*{12mm}on dual spaces to Lie algebras}

For the linear Hamiltonian $H, F \in C_{{\mathcal G}^{*}} = C_u,$ 
\be
H = \pbf{u}^t \pbf{Y}, \qquad  F = \pbf{u}^t \pbf{Z}, 
\ee
The quadratic Poisson brackets (6.3) on ${\mathcal G}^*$ yields:
\be
\{H, F\}_{{\mathcal G}^{*}} \sim \langle \pbf{Y}^. \pbf{u},  {\mathcal O} (\pbf{Z}^. \pbf{u}) \rangle .
\ee
Therefore, the image under the map $\Phi:  {\mathcal O}_{{\mathcal G}^{*}} 
\rightarrow C_{\mathcal M}$ of the Poisson 
bracket $\{H, F\}_{{\mathcal G}^{*}}$ is: 
\be
\Phi (\{H, F\}_{{\mathcal G}^{*}} ) \sim \langle \pbf{Y}^. (\pbf{x} \nabla \pbf{p}), 
{\mathcal O} (\pbf{Z}^. (\pbf{x} \nabla \pbf{p}))\rangle . 
\ee
On the other hand, by formula (7.9), 
\be
\{ \Phi (H), \Phi (F) \}_{\mathcal M} 
\sim  \left(\matrix{-\pbf{Z}^. \pbf{p} \cr \pbf{Z}. \pbf{x} \cr} \right)^t \ \left(\matrix{
B_{xx} & B_{xp} \cr B_{px} & B_{pp} \cr} \right) \ \left(\matrix{-\pbf{Y}^. \pbf{p} \cr 
\pbf{Y}. \pbf{x} \cr} \right), 
\ee
where
\be
B = \ \left(\matrix{B_{xx} & B_{xp} \cr B_{px} & B_{pp} \cr} \right) 
\ee
is the following {\it quadratic} Poisson bracket on $V\oplus V^*$: 
\renewcommand{\theequation}{\arabic{section}.\arabic{equation}{\rm a}}
\setcounter{equation}{5}
\be
(\{\pbf{M}^t \pbf{x}, \pbf{N}^t \pbf{x} \} \sim ) \  \pbf{N}^t B_{xx} (\pbf{M})  \sim \langle \pbf{x} \nabla 
\pbf{M}, {\mathcal O} (\pbf{x} \nabla \pbf{N}) \rangle, 
\ee
\renewcommand{\theequation}{\arabic{section}.\arabic{equation}{\rm b}}
\setcounter{equation}{5}
\be
(\{\pbf{M}^t \pbf{p}, \pbf{N}^t \pbf{x} \} \sim)  \ \pbf{N}^t B_{xp} (\pbf{M})  \sim  - \langle \pbf{M} \nabla 
\pbf{p}, {\mathcal O} (\pbf{x} \nabla \pbf{N}) \rangle , 
\ee
\renewcommand{\theequation}{\arabic{section}.\arabic{equation}{\rm c}}
\setcounter{equation}{5}
\be
(\{\pbf{M}^t \pbf{x}, \pbf{N}^t \pbf{p}\}  \sim) \  \pbf{N}^t B_{px} (\pbf{M})  \sim   - \langle
 \pbf{x} \nabla  \pbf{M},  {\mathcal O} (\pbf{N} \nabla \pbf{p}) \rangle,
\ee
\renewcommand{\theequation}{\arabic{section}.\arabic{equation}{\rm d}}
\setcounter{equation}{5}
\be
(\{\pbf{M}^t \pbf{x}, \pbf{N}^t \pbf{p}\}  \sim) \  \pbf{N}^t B_{pp} (\pbf{M})  \sim   \langle
 \pbf{M} \nabla \pbf{p},  {\mathcal O} (\pbf{N} \nabla \pbf{p}) \rangle,
\ee
The formulae (8.6b) and (8.6c) obviously agree with each other; the matrix $B$  (8.5) is 
thus skewsymmetric.  In the next Section we shall verify that the corresponding quadratic 
Poison bracket on $V \oplus V^*$ satisf\/ies the Jacobi identity and is compatible with the 
symplectic Poisson bracket (7.4).

Let us check now that the map $\Phi: C_{{\mathcal G}^{*}} \rightarrow C_\mu$,
$\Phi (\pbf{u}) = \pbf{x} \nabla \pbf{p}$, is Hamiltonian.  Writing in long hand the 
expression~(8.4) and using formulae~(8.6), we get
\[
\ba{l}
\ds \{ \Phi (H), \Phi(F)\}_{\mathcal M}  \sim  (- \pbf{Z}^. \pbf{p})^t B_{xx} (- \pbf{Y}^. \pbf{p}) 
+ (- \pbf{Z}^. \pbf{p})^t B_{xp} (\pbf{Y}. \pbf{x}) 
\vspace{2mm}\\
\ds \qquad 
+ (\pbf{Z}. \pbf{x})^t B_{px} (- \pbf{Y}^. \pbf{p}) + (\pbf{Z}. \pbf{x})^t B_{pp} (\pbf{Y}. \pbf{x}) 
\vspace{2mm}\\
\ds \qquad \sim \langle \pbf{x} \nabla (\pbf{Y}^. \pbf{p}), {\mathcal O} (\pbf{x} \nabla 
(\pbf{Z}^. \pbf{p})) \rangle + \langle (\pbf{Y}. \pbf{x}) \nabla \pbf{p}, 
 {\mathcal O} (\pbf{x} \nabla (\pbf{Z}^. \pbf{p})) \rangle
\vspace{2mm}\\
\ds \qquad +\langle \pbf{x} \nabla (\pbf{Y}^. \pbf{p}), {\mathcal O} ((\pbf{Z}. \pbf{x}) \nabla \pbf{p}) 
\rangle  + \langle (\pbf{Y}. \pbf{x}) \nabla 
\pbf{p},  {\mathcal O} ((\pbf{Z}.  \pbf{x} ) \nabla \pbf{p}) \rangle
\vspace{2mm}\\
\ds \qquad = \langle \pbf{x} \nabla (\pbf{Y}^. \pbf{p}) + (\pbf{Y}. \pbf{x}) \nabla \pbf{p}, 
 {\mathcal O} (\pbf{x} (\nabla (\pbf{Z}. \pbf{p}) + (\pbf{Z}. \pbf{x}) \nabla \pbf{p}) \rangle
\vspace{2mm}\\
\ds \qquad {\mathop{=}\limits^{\mbox{\scriptsize [by (8.8)]}}} \
 \langle \pbf{Y}^. (\pbf{x} \nabla \pbf{p}),  {\mathcal O} (\pbf{Z}^. (\pbf{x} \nabla \pbf{p})) \rangle \ 
{\mathop{=}\limits^{\mbox{\scriptsize [by (8.3)]}}} \  \Phi (\{H, F\})_{{\mathcal G}^{*}}. 
\ea
\]

\noindent
{\bf Lemma 8.7.}
\renewcommand{\theequation}{\arabic{section}.\arabic{equation}}
\setcounter{equation}{7}
\be 
X^. (v \nabla v^*) = (X. v) \nabla v^* + v \nabla (X^. v^*), \qquad \forall \; v \in V, \ 
v^* \in V^*, \  X \in {\mathcal G}. 
\ee

\noindent
{\bf Proof.}  For any $L \in {\mathcal G}$, we have
\[
\ba{l}
\ds \langle X^. (v \nabla v^*),  L \rangle \sim \langle v \nabla v^*, [L, X] \rangle \sim 
\langle v^*, [L, X]. v \rangle
\vspace{2mm}\\
\ds \qquad = \langle v^*, L. (X. v) - X. (L.v) \rangle \sim \langle (X. v) \nabla v^*, L \rangle 
+ \langle v \nabla (X^. v^*), L \rangle
\vspace{2mm}\\
\ds \qquad =  \langle (X. v) \nabla v^* + v \nabla (X^\cdot v^*), L \rangle. 
\hspace{222.4pt} \mbox{\rule{3mm}{3mm}}
\ea
\]

\noindent
{\bf Remark 8.8.} Originally, Clebsch representation was discovered by Clebsch in 
vector calculus on ${\mathbf R}^n$, $n = 2,3$, without any Lie-algebraic connections.  
The later are the 
results of more recent developments, in 1980's, and are summarized and developed in my 
book~[6].  Since the publications of that book in 1992, there have been~2.5 other developments 
I'm aware of.  First, I found (in~[8] \S~6) {\it quantum} Clebsch representations for the 
linear Poisson brackets on the dual spaces to {\it finite-dimensional} Lie algebras.  Second, 
for the general non-quantal case, the theory of Clebsch representations has been generalized into 
the noncommutative realm in~[10].  Finally, Dr. Morrison in~[12] pp.~500--503 independently 
published a simple version of some of the Clebsch representations results from~[5], but in a 
slightly dif\/ferent notation.  

\medskip

We conclude this Section by writing down explicit quadratic Poisson bracket formulae on 
$V \oplus V^*$ for the case of f\/inite dimensions.  In the notation (6.35)--(6.37), let 
$\{\ell_\alpha\}$ be a basis in $V$, $\{\ell^\alpha\}$ the dual basis in $V^*$, and 
\setcounter{equation}{8}
\be
\chi (e_i) (\ell_\alpha) = \sum_\gamma \chi _{i \alpha}^\gamma \ell_\gamma
\ee
be the action formulae for the representation $\chi: {\mathcal G} \rightarrow  \mbox{End}\,(V)$.  
Then formuale (6.6) yieid
\renewcommand{\theequation}{\arabic{section}.\arabic{equation}{\rm a}}
\setcounter{equation}{9}
\be
\{x^\alpha, x^\beta\} = \sum_{st\mu\nu} r^{st} \chi^\alpha_{s \mu} \chi^\beta_{t \nu} 
x^\mu x^\nu, 
\ee
\renewcommand{\theequation}{\arabic{section}.\arabic{equation}{\rm b}}
\setcounter{equation}{9}
\be
\{p_\alpha, x^\beta\} = - \sum_{st\mu\nu} r^{st} \chi^\mu_{s \alpha} \chi^\beta_{t \nu} 
p_\mu x^\nu, 
\ee
\renewcommand{\theequation}{\arabic{section}.\arabic{equation}{\rm c}}
\setcounter{equation}{9}
\be
\{p_\alpha, p_\beta\} = \sum_{st\mu\nu} r^{st} \chi^\mu_{s \alpha} \chi^\nu_{t \beta} 
p_\mu p_\nu. 
\ee
The map $\Phi(\pbf{u}) = \pbf{x} \nabla \pbf{p}$ takes the form
\renewcommand{\theequation}{\arabic{section}.\arabic{equation}}
\setcounter{equation}{10}
\be
\Phi (u_s) = \sum_{\alpha \beta} \chi^\beta_{s \alpha} x^\alpha p_\beta, 
\ee
and it is a Hamiltonian map between the quadratic Poisson bracket (6.37) on ${\mathcal G}^*$,
\be
\{u_i, u_j\} = \sum_{st \kappa \ell} r^{st} c^\kappa_{is} c^\ell_{jt} u_\kappa u_\ell, 
\ee
and the quadratic Poisson brackets (8.10) on $V \oplus V^*$.

From the results in the next Section it follows that the quadratic Poisson brackets (8.10) 
on $V \oplus V^*$ satisfy the Jacobi identity and are compatible with the symplectic Poisson 
bracket
\be
\{x^\alpha, p_\beta\} = \delta^\alpha_\beta, \qquad  \{x^\alpha, x^\beta\} = \{p_\alpha, p_\beta\} = 0. 
\ee
Finite-dimensional formulae (8.10) can be found in Zakrzewski's paper~[15].

\setcounter{equation}{0}


\section{Properties of the quadratic Poisson brackets on 
{\mathversion{bold}$V \oplus V^*$}}

In this Section we prove that:   1) the quadratic Poisson bracket (8.6) on
$\mbox{Fun}\, (V\oplus V^*$), 
induced by an ${\mathcal O}$-operator ${\mathcal O}: {\mathcal G}^* \rightarrow {\mathcal G}$ 
and a representation $\chi: {\mathcal G} \rightarrow \mbox{Dif\/f}\, (V)$, is legitimate, i.e., 
it satisf\/ies the Jacobi identity;  2)  this quadratic 
Poisson bracket is compatible with the symplectic one;  3) 
 the natural action of ${\mathcal G}$ on 
$\mbox{Fun}\,(V \oplus V^*)$ satisf\/ies the inf\/initesimal Hamiltonian action criterion (5.34) for this 
quadratic Poisson bracket.

\medskip

\noindent
{\bf Proposition 9.1.} {\it  The quadratic Poisson brackets (8.6) on 
$\mbox{\rm Fun} \, (V \oplus V^*)$ satisfy the Jacobi identity. }

\medskip

\noindent
{\bf Proof.}  By the main Theorem of the Hamiltonian formalism, we have to 
verify that 
\setcounter{equation}{1}
\be \{ \{H, F\}, G \} + \mbox{c.p.}  \sim  0
\ee
for all Hamiltonians $H$, $F$, $G$ {\it linear} in the $x$'s 
and the $p$'s.  We break the verif\/ication 
procedure into 4 cases indexed by the number of the $p$'s involved in $H$, $F$, $G$: 
zero, one, two, or three.

\noindent
{\bfseries \itshape Case zero:} 
\be
H = \pbf{X}^t \pbf{x}, \qquad F = \pbf{Y}^t \pbf{x}, \qquad G = \pbf{Z}^t \pbf{x}, 
\ee
where, as understood throughout this paper, $\pbf{X}$, $\pbf{Y}$, $\pbf{Z}$ 
are arbitrary vectors with 
entries in $R$ (or $\tilde R \supset R$).  By formula (8.6a), 
\be
\{H, F\}  \sim \langle \pbf{x} \nabla \pbf{X}, {\mathcal O} (\pbf{x} \nabla \pbf{Y}) \rangle.
\ee
Denoting temporarily
\be
\bar X = \pbf{x} \nabla \pbf{X}, \qquad
 \bar{\!\bar X} = {\mathcal O} (\bar X), 
\ee
we get from formula (9.4) and the relations
\be
\langle  \pbf{x} \nabla \pbf{X}, \bar{\bar Y}\rangle \sim \langle \pbf{X}, 
\bar{\bar Y}.\pbf{x} \rangle \sim  - \langle \bar{\bar Y}{}^{\cdot} \pbf{X}, \pbf{x} \rangle
\ee
that 
\be
{\delta \over \delta \pbf{x}} ( \{H, F \})  \sim - \bar{\bar Y}{}^{\cdot} 
\pbf{X} + \bar{\!\bar X}{}^{\cdot} \pbf{Y}. 
\ee
Therefore, by formula (8.6a) again,
\be
\{ \{H, F\}, G\}  \sim  \pbf{Z}^t B_{xx} (\bar{\!\bar X}{}^{\cdot} 
\pbf{Y} - \bar{\bar Y}{}^{\cdot} \pbf{X})  \sim \langle \pbf{x} \nabla 
(\bar{\!\bar X}{}^{\cdot} \pbf{Y} - \bar{\bar Y}{}^{\cdot} \pbf{X}), 
\bar{\!\bar  Z}\rangle,
\ee
so that
\[
\ba{l}
\ds \{\{ H, F \}, G \}+ \mbox{c.p.}  \sim ( \langle \pbf{x} \nabla (\bar{\!\bar X}{}^{\cdot} \pbf{Y}),
 \bar{\!\bar Z}\rangle  + \mbox{c.p.} ) - ( \langle 
\pbf{x} \nabla (\bar{\bar  Y}{}^{\cdot} \pbf{X}), \bar{\!\bar Z}\rangle + \mbox{c.p.} )
\vspace{2mm}\\
\ds \qquad = ( \langle \pbf{x} \nabla (\bar{\!\bar X}{}^{\cdot} \pbf{Y}), \bar{\!\bar Z }
\rangle + \mbox{c.p.} ) - ( \langle \pbf{x} \nabla (\bar{\!\bar Z}{}^{\cdot} \pbf{Y}), 
\bar{\!\bar  X}\rangle + \mbox{c.p.} )
\vspace{2mm}\\
\ds \qquad = ( \langle \pbf{x} \nabla (\bar{\!\bar X}{}^{\cdot} \pbf{Y}), \bar{\!\bar  Z} 
\rangle - \langle \pbf{x} \nabla (\bar{\!\bar Z}{}^{\cdot} \pbf{Y}), \bar{\!\bar X}\rangle ) 
+ \mbox{c.p.}
\vspace{2mm}\\
\ds \qquad \sim \langle - \bar{\!\bar Z}{}^{\cdot} (\bar{\!\bar X}{}^{\cdot} \pbf{Y}) + 
\bar{\!\bar X}{}^{\cdot}  (\bar{\!\bar Z}{}^{\cdot}
 \pbf{Y}),  \pbf{x} \rangle + \mbox{c.p.} =  \langle [\bar{\!\bar  X}, 
\bar{\!\bar Z}]^\cdot \pbf{Y}, \pbf{x} \rangle + \mbox{c.p.}
\vspace{2mm}\\
\ds \qquad \sim \langle \pbf{x} \nabla \pbf{Y}, [\bar{\!\bar Z}, \bar{\!\bar X}] ) 
\rangle + \mbox{c.p.} = \langle \bar{\bar  Y}, [{\mathcal O} (\bar{\!\bar Z}), {\mathcal O} 
(\bar{\!\bar  X})]\rangle + \mbox{c.p.}  \sim  0 
\ea
\]
by formula (2.6); 

\newpage

\noindent
{\bfseries \itshape Case one:}
\be
H = \pbf{X}^t \pbf{x}, \qquad F = \pbf{Y}^t \pbf{x}, \qquad G = \pbf{Z}^t \pbf{p}. 
\ee
By formulae (9.7) and (8.6c), 
\be
\{ \{H, F\}, G \} \sim \pbf{Z}^t B_{px} ( \bar{\!\bar X}{}^{\cdot} 
\pbf{Y} - \bar{\bar {\pbf{Y}}}{}^{\cdot} \pbf{X} ) \sim  - \langle  \pbf{x} 
\nabla (\bar{\!\bar X}{}^{\cdot} \pbf{Y} - \bar{\bar Y}{}^{\cdot} \pbf{X}),
 {\mathcal O} (\pbf{Z} \nabla \pbf{p} ) \rangle  .
\ee
On the other hand, 
\[
\{F, G \}  \sim  \pbf{Z}^t B_{px} (\pbf{Y}) \sim - \langle  \pbf{x} \nabla \pbf{Y}, {\mathcal O} 
(\pbf{Z} \nabla \pbf{p}) \rangle  \sim 
 \langle  {\mathcal O} (\pbf{Z} \nabla \pbf{p})^\cdot \pbf{Y}, \pbf{x} \rangle 
 \sim  \langle  \pbf{p}, {\bar{\bar Y}}.\pbf{Z} \rangle.
\]
Thus, in the notation
\[
\underline Z = \pbf{Z} \nabla \pbf{p}, \qquad \underline{\underline Z}  = 
{\mathcal O} (\underline  Z), 
\]
we have
\be
{\delta \over \delta \pbf{x}} (\{ F, G \}) = 
\underline{\underline Z}{}^\cdot
\pbf{Y}, \qquad  {\delta \over \delta \pbf{p}} 
(\{ F, G \} ) = {\bar{\bar Y}} .\pbf{Z}. 
\ee
Therefore, by formulae (8.6a,b), 
\[
\{ \{F, G\}, H \} \sim \pbf{X}^t[B_{xx} (\underline{\underline Z}{}^\cdot 
\pbf{Y}) + B_{xp} ({\bar{\bar Y}} .\pbf{Z}) ] 
\]
\be 
\qquad \sim  \langle  \pbf{x} \nabla (\underline{\underline Z} \pbf{Y}), 
\bar{\!\bar X} \rangle  - \langle 
 ({\bar{\bar Y}} .\pbf{Z}) \nabla \pbf{p}, 
\bar{\!\bar X} \rangle .
\ee
Interchanging $X$ and $Y$ in formula (9.12), we obtain 
\[
\{ \{G, H \}, F\}  \sim  - \{ \{ H, G \}, F \} 
\]
\be
\qquad \sim  - \langle  \pbf{x} \nabla (\underline{\underline Z}{}^\cdot \pbf{X}), 
\bar{\bar Y} \rangle  + \langle  (\bar{\!\bar X}{}^{\cdot} \pbf{Z}) \nabla \pbf{p}, 
\bar{\bar Y} \rangle  .
\ee
The 2$^{nd}$ summands in the expressions (9.12) and (9.13) combine into
\[
\langle  \pbf{p}, - {\bar{\!\bar  X}}.({\bar{\bar Y}}.\pbf{Z}) + 
{\bar{\bar  Y}}.({\bar{\!\bar X}}.\pbf{Z}) \rangle   =  \langle  \pbf{p}, 
[\bar{\bar Y}, \bar{\!\bar  X}]. \pbf{Z} \rangle  
\]
\[
\qquad \sim  \langle  \pbf{Z} \nabla \pbf{p}, [ {\mathcal O} (\overline Y),
 {\mathcal O} (\overline X)] \rangle  
\ {\mathop{=}\limits^{\mbox{\scriptsize [by (2.10)]}}}  \ 
\langle  \underline Z, {\mathcal O} (\bar{\bar Y}{}^{\cdot} \bar
 X - \bar{\!\bar X}{}^{\cdot} \bar Y) \rangle  
\]
\be
\qquad \sim - \langle  \bar{\bar Y}{}^{\cdot} \bar X - 
\bar{\!\bar X}{}^{\cdot}  Y, 
\underline{\underline Z}\rangle ,
\ee
while the 1$^{st}$ summands in the expressions (9.12) and (9.13) combine into 
\[
\langle  \underline{\underline Z}{}^\cdot 
\pbf{Y}, \underline{\underline X} . \pbf{x} \rangle  - 
\langle  \pbf{Z}^\cdot \pbf{x}, 
\underline{\underline Y} . \pbf{x} \rangle   \sim - 
\langle  \pbf{Y}, \underline{\underline Z}  . 
(\underline{\underline X} . \pbf{x}) \rangle  
+ \langle  \pbf{X},  \underline{\underline Z}  . 
(\underline{\underline Y}   . \pbf{x}) \rangle 
\]
\be
\qquad \sim  - \langle  ({\bar{\!\bar X}}.\pbf{x}) \nabla \pbf{Y} + (
{\bar{\bar Y}}.\pbf{x}) \nabla \pbf{X}, \bar{\!\bar Z} \rangle  . 
\ee
We see that the sum total of the expressions (9.10), (9.14), and (9.15) is $\sim 0$ provided
\[
- \pbf{x} \nabla (\bar{\!\bar X}{}^{\cdot} \pbf{Y} - 
\bar{\bar Y}{}^{\cdot} \pbf{X}) - \bar{\bar Y}{}^{\cdot} 
\bar X + \bar{\!\bar X}{}^{\cdot} \bar Y - 
(\underline{\underline X}. \pbf{x}) \nabla \pbf{Y} +
 (\underline{\underline Y}. \pbf{x}) \nabla \pbf{X} = 0, 
\]
and this is so by formulae (8.8) and (9.5); 

\noindent
{\bfseries \itshape Cases two and three}
 follow from cases one and zero, respectively, once we notice 
that the quadratic Poisson bracket formulae  (8.6) allow the symmetry
\be
\pbf{x} \mapsto \pbf{p}, \qquad  \pbf{p} \mapsto \pbf{x}, \qquad
 \chi \mapsto \chi^d,
\ee
where $\chi^d: {\mathcal G} \mapsto \mbox{Dif\/f}\, (V^*)$ is the dual representation,
\be
\chi^d (X) = - \chi (X)^\dagger, \qquad   \forall \; X \in {\mathcal G}. 
\ee
This symmetry  becomes obvious if we use the relation 
\be
a \nabla b = - b {\mathop{\nabla}\limits_d} a
\ee
where ${\mathop{\nabla}\limits_d}$ product is taken w.r.t. 
the dual representation $\chi^d$:
\be
\langle  a \nabla b, X \rangle  \sim  \langle  b, X . a \rangle   \sim 
 \langle  - X^\cdot b, a \rangle  =  \langle a,  - X^\cdot b \rangle  
\sim  \langle  - b {\mathop{\nabla}\limits_d} a, X\rangle  , \qquad  \forall \;
X \in {\mathcal G}.
\ee
Formulae (8.6a) and (8.6d) are interchanged under the symmetry (9.16), as are formulae 
(8.6b) and (8.6c). \hfill \rule{3mm}{3mm}

\medskip

\noindent
{\bf Proposition 9.20.} {\it The quadratic Poisson bracket (8.6) on $V \oplus V^*$ and 
the symplectic bracket (7.4) are compatible.}

\medskip

\noindent
{\bf Proof.}  We have to verify that 
\setcounter{equation}{20}
\be
( \{\{ H, F \}_1, G\}_2 + \{ \{ H, F \}_2, G \}_1) + \mbox{c.p.} \sim 0 
\ee
for all $H$, $F$, $G$ linear in $\pbf{x}$, $\pbf{p}$, 
where $\{ \ , \}_1$ denotes the symplectic Poisson 
bracket and $\{ \ ,  \}_2$ denotes the quadratic one.  This is obviously true  for $H$, $F$, $G$ 
linear in $\pbf{x}$; and also for $H$, $F$, $G$ linear in $\pbf{p}$.  

\noindent
{\bfseries \itshape Case one:} 
\be
H = \pbf{X}^t \pbf{x}, \qquad F = \pbf{Y}^t \pbf{x}, \qquad
 G =  \pbf{Z}^t \pbf{p}. 
\ee
We have
\[
\{H, F\}_1 = 0, \qquad \{F, G\}_1 \sim - \pbf{Y}^t \pbf{Z}, \qquad \{G, H \}_1 \sim \pbf{Z}^t \pbf{X}. 
\]
Thus,
\[
\{ \{ ( \ \cdot \ ), (\cdot \cdot)\}_1, ( \cdot \cdot \cdot)\}_2 = 0. 
\]
Now, by formula (9.7),
\be
{\delta \{H, F \} \over \delta \pbf{x}} = \bar{\!\bar  X}{}^{\cdot} 
\pbf{Y} - \bar{\bar Y}{}^{\cdot} \pbf{X}, 
\ee
so that, by formulae (7.4), 
\be
\{ \{H, F\}_2, G\}_1 \sim - \pbf{Z}^t (\bar{\!\bar X}{}^{\cdot} 
\pbf{Y} - \bar{\bar Y}{}^{\cdot} \pbf{X}) \sim \pbf{Y}^t ({\bar{\!\bar X}}.\pbf{Z}) 
- \pbf{X}^t ({\bar{\bar Y}}.\pbf{Z}).
\ee
Next, by formula (9.11),
\[
{\delta \{ F, G\}_2 \over \delta \pbf{p}} = {\bar{\bar Y}}.\pbf{Z},
\]
so that
\renewcommand{\theequation}{\arabic{section}.\arabic{equation}{\rm a}}
\setcounter{equation}{24}
\be
\{ \{F, G\}_2, H \}_1 \sim \pbf{X}^t ({\bar{\bar Y}}.\pbf{Z}).
\ee
Interchanging $X$ and $Y$ in the formula (9.25a), we obtain 
\renewcommand{\theequation}{\arabic{section}.\arabic{equation}{\rm b}}
\setcounter{equation}{24}
\be
\{\{G, H\}_2, F\}_1 \sim - \{ \{H, G\}_2, F\}_1 \sim - \pbf{Y}^t 
({\bar{\!\bar X}}.\pbf{Z}).
\ee
Adding up the expression (9.24)--(9.26) we get zero.

\noindent
{\bfseries \itshape Case two}  follows from the already established case one (9.22) by the 
application of the symmetry (9.16) accompanied by the replacement of the symplectic matrix 
$b$ (7.3) by $-b$, -- which doesn't af\/fect the validity of the case one arguments.
\hfill \rule{3mm}{3mm}

\medskip

\noindent
{\bf Remark 9.26.} The Proof of Proposition 9.1 could be reduced to only  the 
case zero upon noticing that formulae (8.6b--d) are particular instances of the basic formula 
(8.6a) applied to the representation $\chi^{\mbox{\scriptsize new}} = 
\chi \oplus \chi^d $ on $V^{\mbox{\scriptsize new}} = V \oplus V^*$. 

\medskip

\noindent
{\bf Proposition 9.27.} {\it The natural (anti) action of ${\mathcal G}$ on $C_{\mathcal M} = 
\mbox{\rm Fun} \, (V \oplus V^*)$, 
\renewcommand{\theequation}{\arabic{section}.\arabic{equation}}
\setcounter{equation}{27}
\be
X^\wedge (\pbf{x}) = X. \pbf{x}, \qquad  X^\wedge (\pbf{p}) = X^\cdot \pbf{p}, 
\ee
satisfies the criterion (5.34) of infinitesimal Hamiltonian action w.r.t. the quadratic 
Poisson bracket (8.6) on $V \oplus V^*$.}

\medskip

\noindent
{\bf Proof.} We shall check formula (5.34) for the linear Hamiltonians 
$H$, $F$.  We break this check into 3~cases,
 depending upon how many $p$'s are present among $H$ and $F$. 

\noindent
{\bfseries \itshape Case zero:}
\renewcommand{\theequation}{\arabic{section}.\arabic{equation}}
\setcounter{equation}{28} 
\be
H = \pbf{Y}^t \pbf{x}, \qquad  F = \pbf{Z}^t \pbf{x}. 
\ee
We have, by formula (8.6a):
\be
\ba{l}
X^\wedge (\{H, F\}) \sim X^\wedge (\langle  \pbf{x} \nabla \pbf{Y}, 
{\mathcal O} (\pbf{x} \nabla \pbf{Z}) \rangle  ) 
\vspace{2mm}\\
\qquad = \langle  (X . \pbf{x}) \nabla \pbf{Y}, {\mathcal O} (\pbf{x} \nabla \pbf{Z}) \rangle  
+ \langle  \pbf{x} \nabla \pbf{Y}, {\mathcal O} (( X . \pbf{x}) 
\nabla \pbf{Z}) \rangle  ,
\ea
\ee
\be
- \{X^\wedge (H), F \} \sim - \{ \pbf{Y}^t (X . \pbf{x}), \pbf{Z}^t \pbf{x}\} \sim \langle  \pbf{x} \nabla 
(X^\cdot \pbf{Y}), {\mathcal O} (\pbf{x} \nabla \pbf{Z}) \rangle , 
\ee
\be
-\{H, X^\wedge (F)\} = - \{ \pbf{Y}^t \pbf{x}, \pbf{Z}^t (X. \pbf{x})\} \sim 
 \langle  \pbf{x} \nabla \pbf{Y}, {\mathcal O} ((X. \pbf{x}) \nabla \pbf{Z}) \rangle  . 
\ee
Adding up the expressions (9.30)--(9.32) and using formula (8.8), we f\/ind
\be
\ba{l}
X^\wedge (\{H, F\}) - \{ X^\wedge (H), F \} - \{H, X^\wedge (F) \}
\vspace{2mm}\\
\qquad \sim  \langle  X^\cdot (\pbf{x} \nabla \pbf{Y}), {\mathcal O} (\pbf{x} \nabla \pbf{Z}) \rangle 
 + \langle  \pbf{x} \nabla \pbf{Y}, {\mathcal O} 
(X^\cdot (\pbf{x} \nabla \pbf{Z})) \rangle  . 
\ea
\ee
On the other hand, formulae (9.28) imply that 
\renewcommand{\theequation}{\arabic{section}.\arabic{equation}{\rm a}}
\setcounter{equation}{33}
\be
X^\wedge (\pbf{Y}^t \pbf{x}) = \pbf{Y}^t (X. \pbf{x}) \sim
 \langle  \pbf{x} \nabla \pbf{Y}, X \rangle , 
\ee
\renewcommand{\theequation}{\arabic{section}.\arabic{equation}{\rm b}}
\setcounter{equation}{33}
\be
X^\wedge (\pbf{Y}^t \pbf{p}) = \pbf{Y}^t (X^\cdot \pbf{p}) \sim 
 \langle  - \pbf{Y} \nabla \pbf{p}, X \rangle  , 
\ee
so that 
\renewcommand{\theequation}{\arabic{section}.\arabic{equation}}
\setcounter{equation}{34}
\be
(\pbf{Y}^t \pbf{x})^\sim = \pbf{x} \nabla \pbf{Y}, \qquad
 (\pbf{Y}^t \pbf{p})^\sim = - \pbf{Y} \nabla \pbf{p}. 
\ee
Therefore,
\be
\ba{l}
\langle  [H^\sim, F^\sim], X) = \langle [\pbf{x} \nabla \pbf{Y}, \pbf{x} \nabla \pbf{Z}],  X \rangle 
\vspace{2mm}\\
\qquad {\mathop{=}\limits^{\mbox{\scriptsize [by (2.11) ]}}} \
 \langle  {\mathcal O}( \pbf{x} \nabla \pbf{Y}). (\pbf{x} \nabla \pbf{Z}) 
- {\mathcal O} (\pbf{x} \nabla \pbf{Z}). (\pbf{x} \nabla 
\pbf{Y}), X \rangle 
\vspace{2mm}\\
\qquad \sim  \langle  - X^. (\pbf{x} \nabla \pbf{Z}),  {\mathcal O} (\pbf{x} \nabla \pbf{Y}) \rangle  
+ \langle  X^\cdot (\pbf{x} \nabla Y), 
{\mathcal O} (\pbf{x} \nabla \pbf{Z}) \rangle  , 
\ea
\ee
and this expression is $\sim$ to (9.33) since ${\mathcal O}$ is skewsymmetric.  

\noindent
{\bf Case one:} 
\be
H = \pbf{Y}^t \pbf{p}, \qquad   F = \pbf{Z}^t \pbf{x}.
\ee
We have, by formula (8.6b): 
\be
\ba{l}
X^\wedge (\{H, F\}) \sim X^\wedge (- \langle  \pbf{Y} \nabla \pbf{p}, {\mathcal O} (\pbf{x} 
\nabla \pbf{Z})\rangle  ) 
\vspace{2mm}\\
\qquad = - \langle  \pbf{Y} \nabla (X^\cdot \pbf{p}), {\mathcal O} (\pbf{x} \nabla \pbf{Z}) \rangle  
- \langle  \pbf{Y} \nabla \pbf{p}, {\mathcal O} (
(X. \pbf{x}) \nabla \pbf{Z}) \rangle  , 
\ea
\ee
\be
\ba{l}
- \{ X^\wedge (H), F\} = - \{\pbf{Y}^t (X^\cdot \pbf{p}), \pbf{Z}^t \pbf{x} \}
\vspace{2mm}\\
\qquad  \sim \{ (X. \pbf{Y})^t 
\pbf{p}, \pbf{Z}^t \pbf{x} \} 
\sim - \langle  (X. \pbf{Y}) \nabla \pbf{p}, {\mathcal O} (\pbf{x} \nabla \pbf{Z}) \rangle  , 
\ea
\ee
\be
\ba{l}
- \{ H, X^\wedge (F) \} = - \{ \pbf{Y}^t \pbf{p}, \pbf{Z}^t (X. \pbf{x}) \} 
\vspace{2mm}\\
\qquad \sim \{ \pbf{Y}^t \pbf{p}, 
(X^. \pbf{Z})^t \pbf{x} \} 
\sim - \langle  \pbf{Y} \nabla \pbf{p}, {\mathcal O} (\pbf{x} \nabla (X^\cdot \pbf{Z})) \rangle  . 
\ea
\ee
Adding the expressions (9.38)--(9.40) up and using formula (8.8), we f\/ind
\be
\ba{l}
X^\wedge (\{H, F \}) - \{X^\wedge (H), F\} - \{H, X^\wedge (F) \} 
\vspace{2mm}\\
\qquad \sim - \langle  \pbf{Y} \nabla \pbf{p}, {\mathcal O} (X^\cdot (\pbf{x} \nabla \pbf{Z})) \rangle 
 - \langle  X^\cdot (\pbf{Y} \nabla 
\pbf{p}), {\mathcal O} (\pbf{x} \nabla \pbf{Z}) \rangle . 
\ea
\ee
On the other hand, by formulae (9.35) and (2.11), 
\be
\ba{l}
\langle  [H^\sim, F^\sim ], X \rangle  = - \langle  [\pbf{Y} \nabla \pbf{p}, \pbf{x} \nabla \pbf{Z}], 
X\rangle  = - \langle  {\mathcal O} (\pbf{Y} 
\nabla \pbf{p}). (\pbf{x} \nabla \pbf{Z}), X \rangle  
\vspace{2mm}\\
\qquad + \langle  {\mathcal O} (\pbf{x} \nabla \pbf{Z}). (\pbf{Y} \nabla \pbf{p}), X \rangle   \sim 
 \langle   X^\cdot (\pbf{x} \nabla \pbf{Z}), 
{\mathcal O} (\pbf{Y} \nabla \pbf{p}) \rangle 
- \langle  X^\cdot (\pbf{Y} \nabla \pbf{p}), {\mathcal O} (\pbf{x} \nabla \pbf{Z}) \rangle , 
\ea\hspace{-6.25pt}
\ee
and this is $\sim$ to the expression (9.41) because ${\mathcal O}$ is skewsymmetric.

\noindent
{\bfseries \itshape Case two}
 follows from the case zero and the symmetry property (9.16). \hfill \rule{3mm}{3mm}

\renewcommand{\theequation}{{\rm A1}.\arabic{equation}}
\setcounter{equation}{0}

\section*{Appendix A1.  Crossed Lie algebras}

Let $G$ be a Hamilton-Lie group, i.e., a Lie group with a multiplicative Poisson bracket 
on it.  The inf\/initesimal version of this object is a Lie bialgebra, i.e., a Lie bracket [~,~] 
on ${\mathcal G}^*$ whose dual, considered as a map 
$\varphi: {\mathcal G} \rightarrow \wedge^2 {\mathcal G}$, is 1-cocycle 
on ${\mathcal G}$.  Drinfel'd noticed in~[2] that the 1-cocycle condition 
can be reformated in such a 
way as to make the self-duality of the notion of Lie bialgebra explicit, as follows.  Since 
${\mathcal G}$ is a Lie algebra, it acts on its dual space, ${\mathcal G}^*$.  
Consider the following skew multiplication 
on the space ${\mathcal G} + {\mathcal G}^*$:
\be
\left[\left(\matrix{X \cr u \cr} \right), \left(\matrix{Y \cr v \cr} \right) \right] = 
\left(\matrix{[X, Y] + u^\cdot Y - v^\cdot X \cr [u, v] + X^\cdot v - Y^\cdot u \cr} \right), 
\qquad X, Y \in {\mathcal G}, \quad  u, v \in {\mathcal G}^*.
\ee
The Drinfel'd observation mentioned above is that the bracket (A1.1) satisf\/ies the Jacobi 
identity if\/f $\varphi: {\mathcal G} \rightarrow \wedge^2 {\mathcal G}$ is a 1-cocycle  
on ${\mathcal G}$ (see~[1] p.~27).  
Let us f\/ind the form of the condition, originally written down by Drinfel'd in~[2] for the 
f\/inite-dimensional case, equivalent to formula (A1.1) def\/ining a Lie algebra.  

For the ${\mathcal G}$-component of the expression
\be
\left[ \bigg[ \left(\matrix{X \cr u \cr} \right), \left(\matrix{Y \cr v \cr} \right) \right], 
\left(\matrix{Z \cr w \cr} \right) \bigg] + \mbox{c.p.} = \left[ \left(\matrix{[X, Y] + u^\cdot 
Y - v^\cdot X \cr [u, v] + X^\cdot v - Y^\cdot u \cr} \right), \left(\matrix{Z \cr w \cr} 
\right) \right] + \mbox{c.p.} 
\ee
we get:
\renewcommand{\theequation}{{\rm A1}.\arabic{equation}{\rm a}}
\setcounter{equation}{2}
\be
([[X, Y], Z] + c.p.) + ([u^\cdot Y - v^\cdot X, Z] + c.p.) + ([u, v]^\cdot Z + \mbox{c.p.})  
\ee
\renewcommand{\theequation}{{\rm A1}.\arabic{equation}{\rm b}}
\setcounter{equation}{2}
\be
\qquad 
+ ([X^\cdot v - Y^\cdot u)^\cdot Z - w^\cdot [X, Y] + \mbox{c.p.}) - (w^\cdot (u^\cdot Y - v^\cdot 
X) + \mbox{c.p.}).
\ee
The f\/irst summand vanishes by the Jacobi identity in ${\mathcal G}$.  
The 3$^{rd}$ and 5$^{th}$ summands  combine into 
\[
[w, u]^\cdot Y - w^\cdot (u^\cdot Y) + u^\cdot (w^\cdot Y) + \mbox{c.p.}
\]
and this vanishes because the action of ${\mathcal G}^*$ on ${\mathcal G}$ 
is a representation of the Lie 
algebra structure on ${\mathcal G}^*$.  This leaves us with 2$^{nd}$ 
and 4$^{th}$ summands, which combine into
\renewcommand{\theequation}{{\rm A1}.\arabic{equation}}
\setcounter{equation}{3}
\be
[u^\cdot Y, Z] - [u^\cdot Z, Y] + (Z^\cdot u)^\cdot Y - (Y^\cdot u)^\cdot Z - u^\cdot 
[Y, Z] = 0. 
\ee
This is an equation in ${\mathcal G}$.  Let us calculate the value 
$\langle  v, LHS \rangle $ for an arbitrary 
element $v \in {\mathcal G}^*$.  Term-by-term, we f\/ind:
\renewcommand{\theequation}{{\rm A1}.\arabic{equation}{\rm a}}
\setcounter{equation}{4}
\be
(v, [u^\cdot Y, Z] \rangle   \sim \ \langle  Z^\cdot v, u^\cdot Y\rangle  ,
\ee
\renewcommand{\theequation}{{\rm A1}.\arabic{equation}{\rm b}}
\setcounter{equation}{4}
\be 
\langle v, - [u^\cdot Z, Y] \rangle   \sim - \langle  Y^\cdot v, u^\cdot Z \rangle ,
\ee
\renewcommand{\theequation}{{\rm A1}.\arabic{equation}{\rm c}}
\setcounter{equation}{4}
\be
\langle v, (Z^\cdot u)^\cdot Y \rangle   \sim  \langle  [v, Z^\cdot u], Y \rangle  
 \sim - \langle  Z^\cdot u, v^\cdot Y \rangle , 
\ee
\renewcommand{\theequation}{{\rm A1}.\arabic{equation}{\rm d}}
\setcounter{equation}{4}
\be
\langle v, - (Y^\cdot u)^\cdot Z \rangle   \sim  \langle  [Y^\cdot u, v], Z \rangle  
 \sim - \langle  Y^\cdot u, v^\cdot Z \rangle , 
\ee
\renewcommand{\theequation}{{\rm A1}.\arabic{equation}{\rm e}}
\setcounter{equation}{4}
\be
\langle v, -u^\cdot [Y, Z] \rangle   \sim \ \langle  [u, v], [Y, Z] \rangle . 
\ee
Adding the expressions (A1.5) up, we arrive at the following quadrilinear relation equivalent 
to the trilinear equation (A1.4):  
\renewcommand{\theequation}{{\rm A1}.\arabic{equation}}
\setcounter{equation}{5}
\be
\ba{l}
\langle [u, v], [Y, Z] \rangle   \sim  
 \langle  Z^\cdot v, u^\cdot Y\rangle  + \langle Y^\cdot u, v^\cdot Z\rangle  - 
\langle Y^\cdot v, u^\cdot Z\rangle  - \langle Z^\cdot u, v^\cdot Y\rangle,
\vspace{2mm}\\
\forall \; Y, Z \in {\mathcal G}, \qquad  \forall \; u, v \in {\mathcal G}^*  .
\ea
\ee
In this form, the symmetry between ${\mathcal G}$ and ${\mathcal G}^*$ is apparent.  
We don't have to analyze 
the ${\mathcal G}^*$-component of the double-commutator (A1.2).

\medskip

\noindent
{\bf Remark A1.7.} In f\/inite dimenions, the commutator (A1.1) leaves the natural 
scalar product on ${\mathcal G} + {\mathcal G}^*$:  
\setcounter{equation}{7}
\be
\left({X \choose u} , {Y \choose v}\right) = \langle u, Y\rangle  + \langle  v, X \rangle,
\ee
ad-invariant:
\be
\left( \left[ {X \choose u}, {Y \choose v} \right], {Z \choose w}\right)  = 
\left( {X \choose u}, 
\left[{Y \choose v}, {Z \choose w} \right]\right ).
\ee
This ad-invariance is still true in inf\/inite dimensions, provided it is properly understood:
\be
\left(\left[ {X \choose u}, {Y \choose v}\right], {Z \choose w}\right) \sim \left( {X \choose u},
\left[ {Y \choose v}, {Z \choose w} \right] \right).
\ee
Indeed, for the LHS of (A1.10) we get
\[
\ba{l}
\ds \left( \left(\matrix{[X, Y] + u^\cdot Y - v^\cdot X \cr [u, v] + X^\cdot v - Y^\cdot 
u \cr} \right), {Z \choose w} \right) 
\vspace{3mm}\\
\qquad = \langle  [u, v] + X^\cdot u - Y^\cdot u, Z \rangle  +  \langle  w, [X, Y] +
 u^\cdot Y - v^\cdot X \rangle  
\vspace{2mm}\\
\qquad \sim  \langle  u, v^\cdot Z \rangle  - \langle  Z^\cdot v, X\rangle  + 
\langle  u, [Y, Z]\rangle  + 
 \langle  Y^\cdot w, X\rangle  - \langle u, w^\cdot Y \rangle  + \langle  [v, w], X \rangle  
\vspace{2mm}\\
\qquad = \langle  u, v^\cdot Z + [Y,Z] - w^\cdot Y \rangle  + \langle  - Z^\cdot v 
+ Y^\cdot w + [v, w], X \rangle  
\vspace{3mm}\\
\ds \qquad = \left( {X \choose u},  \left(\matrix{[Y, Z] + v^\cdot Z - w^\cdot Y \cr [ v, w ]
+  Y^\cdot w - Z^\cdot v \cr} \right) \right), 
\ea
\]
and this is the RHS of (A1.10).

Let us check that when the commutator on ${\mathcal G}^*$ is given by the formula (2.11), 
\be
[u, v] = {\mathcal O} (u)^\cdot v -  {\mathcal O}(v)^\cdot u,
\ee
where ${\mathcal O}: {\mathcal G}^* \rightarrow {\mathcal G}$ is a skewsymmetric 
${\mathcal O}$-operator, then ${\mathcal O} + {\mathcal O}^*$ 
(A1.1) is a Lie algebra.  Let us check the criterion (A1.4).

\medskip

\noindent
{\bf Lemma A1.12.} 
\setcounter{equation}{12}
\be
u^\cdot X = [ {\mathcal O} (u), X] + {\mathcal O} (X^\cdot u), \qquad
 \forall \; u \in {\mathcal G}^*, \  X \in {\mathcal G}. 
\ee

\noindent
{\bf Proof.}  For any $v \in {\mathcal G}^*$, 
\[
\ba{l}
\langle  v, u^\cdot X \rangle   \sim  \langle  [v, u], X \rangle  = 
 \langle  {\mathcal O} (v)^\cdot u - {\mathcal O} (u)^\cdot v, X \rangle  
 \sim  \langle  u, [ X, {\mathcal O} (v)] \rangle  + \langle  v, [{\mathcal O} (u), X] \rangle
\vspace{2mm}\\
\qquad  \sim - \langle  X^\cdot u, {\mathcal O} (v) \rangle 
 + \langle  v, [{\mathcal O} (u), X] \rangle   \sim 
 \langle  v, {\mathcal O} (X^\cdot u) \rangle  + \langle  v, [{\mathcal O} (u), X] \rangle. 
\hspace{60.5pt} \rule{3mm}{3mm}
\ea
\]

Using formula (A1.13), we transform each of the 5 summands in the LHS of the criterion (A1.4):
\be
1) \quad [u^\cdot Y, Z ] = [[{\mathcal O} (u), Y], Z] + [{\mathcal O} (Y^\cdot u), Z ], 
\ee
\be 
2) \quad - [u^\cdot Z, Y] = - [[{\mathcal O} (u), Z], Y] - [{\mathcal O} (Z^\cdot u), Y], 
\ee
\be
3) \quad (Z^\cdot u)^\cdot Y = [{\mathcal O} (Z^\cdot u), Y] + {\mathcal O}(Y^\cdot (Z^\cdot u)), 
\ee
\be
4) \quad - (Y^\cdot u)^\cdot Z = - [{\mathcal O} (Y^\cdot u), Z] - {\mathcal O} (Z^\cdot (Y^\cdot u)), 
\ee
\be
5) \quad - u^\cdot [Y, Z] = - [{\mathcal O} (u), [Y, Z]] - {\mathcal O} ([Y, Z]^\cdot u). 
\ee
The 2$^{nd}$ summands in the expressions (A1.16)--(A1.18) add up to zero; the $1^{st}$ summands in 
the expressions (A1.14), (A1.15), (A1.18) do likewise; the 2$^{nd}$ summand in the expression 
(A1.14$+ \ell$) and the $1^{st}$ summand in the expression (A1.17$-\ell$), $\ell = 0, 1$, 
cancel each other out.
 
We conclude this Section by examining when the symplectic form on ${\mathcal G} + {\mathcal G}^*$: 
\be
\omega \left( {X \choose u} , {Y \choose v} \right) 
= \langle u, Y\rangle  - \langle v, X \rangle
\ee
is a 2-cocycle on the Lie algebra ${\mathcal G} + {\mathcal G}^*$  (A1.1).  We have:
\[
\omega \left( \left[ {X \choose u}, {Y \choose v} \right], {Z \choose w} \right) + \mbox{c.p.} = 
\omega \left( \left(\matrix{[X, Y] + u^\cdot Y - v^\cdot X \cr [u, v] + X^\cdot v - 
Y^\cdot u \cr} \right), {Z \choose w} \right) +\mbox{c.p.}
\]
\[
\qquad = (\langle  [u, v] + X^\cdot v - Y^\cdot u, Z \rangle  + \mbox{c.p.} )
 - ( \langle  w, [X, Y] + u^\cdot Y - v^\cdot X \rangle  + \mbox{c.p.})
\]
\be
\qquad = (\langle  [u, v], Z \rangle  - \langle u, v^\cdot Z\rangle  + \langle v, u^\cdot Z\rangle  ) + 
\mbox{c.p.}
\ee
\be
\qquad + ( \langle  Y^\cdot w, X\rangle  - \langle  X^\cdot w, Y \rangle  - \langle  w, [X, Y] \rangle ) + 
\mbox{c.p.}
\ee
Thus, the symplectic form (A1.19) is a 2-cocycle if\/f 
\be
\langle  [u, v], Z \rangle  - \langle u, v^\cdot Z \rangle  + \langle v, u^\cdot Z \rangle   \sim 0, 
\qquad  \forall \; u, v \in {\mathcal G}^*, \ Z \in {\mathcal G}, 
\ee
\be
\langle Y^\cdot w, X\rangle  - \langle  X^\cdot w, Y\rangle  - \langle  w, [X, Y] \rangle  \sim 0, 
\qquad  \forall  \; X, Y \in {\mathcal G}, \ w \in {\mathcal G}^*. 
\ee
This happens if\/f, respectively,
\be
[u, v] = 0, \qquad  \forall \; u, v \in {\mathcal G}^*, 
\ee
\be
[X, Y] = 0, \qquad 
 \forall X, Y \in {\mathcal G}, 
\ee
which is this side of ``never''.  The next Section provides a remedy of sorts.

\renewcommand{\theequation}{{\rm A2}.\arabic{equation}}
\setcounter{equation}{0}

\section*{Appendix A2.  Symplectic {\mathversion{bold}$r$}-matrices and symplectic doubles}

Let ${\mathcal G}$ be a Lie algebra and $\rho: {\mathcal G} \rightarrow  \mbox{Dif\/f}\,
({\mathcal G}^*)$ a representation, {\it not  necessarily}  
the coadjoint one.  Let ${\mathcal H}: = {\mathcal G} \ltimes {\mathcal G}^* =
 {\mathcal G} \ {\mathop{\ltimes}\limits_\rho} \ {\mathcal G}^*$ be the 
semidirect sum Lie algebra, with the commutator
\be
\left[ {X \choose u}, {Y \choose v} \right] = \left(\matrix{[X, Y] \cr 
\rho (X) (v) - \rho (Y) (u) \cr} \right). 
\qquad X, Y \in {\mathcal G}, \quad  u, v \in {\mathcal G}^*.
\ee
Lew $\omega$ be the symplectic form on ${\mathcal H} = {\mathcal G} \ltimes {\mathcal G}^*$: 
\be
\omega \left( {X \choose u} , {Y \choose v} \right) = \langle u, Y\rangle  - \langle v, X\rangle .
\ee

\noindent
{\bf Proposition A2.3.} {\it The symplectic form $\omega$ (A2.2) is a 2-cocycle on the 
Lie algebra ${\mathcal G} \ {\mathop{\ltimes}\limits_\rho} \ {\mathcal G}^*$ if\/f the dual 
representation $\rho^d:  {\mathcal G} \rightarrow \mbox{\rm Dif\/f}\, ({\mathcal G})$ 
satisf\/ies the property
\setcounter{equation}{3}
\be
\rho^d (X) (Y) - \rho^d (Y) (X) = [X, Y], \qquad  \forall\; X, Y \in {\mathcal G}. 
\ee}

\noindent
{\bf Proof.} We have,
\be
\ba{l}
\ds \omega \left( \left[ {X \choose u}, {Y \choose v} \right],  {Z \choose w} \right) + \mbox{c.p.}
 = \omega 
\left( \left(\matrix{[X, Y] \cr \rho (X) (v) - \rho (Y) (u) \cr} \right), {Z \choose w} \right) + 
\mbox{c.p.}\!
\vspace{3mm}\\
\qquad = (\langle  \rho (X) (v) - \rho (Y) (u), Z\rangle  - \langle  w, [X, Y] \rangle  ) + \mbox{c.p.}
\vspace{2mm}\\
\qquad = (\langle  \rho (Y) (w), X \rangle  - \langle  \rho (X) (w), Y \rangle  
- \langle  w, [X, Y] \rangle  ) + \mbox{c.p.}
\ea
\ee 
Thus, $\omega$ is a 2-cocycle if\/f
\be
\langle  \rho (Y) (w), X\rangle  - \langle  \rho (X) (w), Y \rangle 
 \sim \langle  w, [X, Y]\rangle , \qquad  \forall \; X, Y \in {\mathcal G}, \ 
w \in {\mathcal G}^*. 
\ee
Rewriting the LHS of this relation as
\[
\langle  w, - \rho^d (Y) (X) \rangle  + \langle  w, \rho^d (X) (Y) \rangle , 
\]
we arrive at the equivalent to (A2.6) equation (A2.4). \hfill \rule{3mm}{3mm}

Assuming from now on that a representation $\rho^d: {\mathcal G} \rightarrow 
\mbox{Dif\/f}\, ({\mathcal G})$ satisfying (A2.4) is f\/ixed, denote by 
\be
xy = \rho^d (x) (y) 
\ee
the resulting multiplication ${\mathcal G} \times {\mathcal G} \rightarrow {\mathcal G}$.  
Since $\rho^d$ is a representation, 
\be
\rho^d ([x, y]) = [\rho^d (x), \rho^d (y) ], \qquad  \forall \; x, y \in {\mathcal G}. 
\ee
Applying this operator identity to an element $z \in {\mathcal G}$, we get
\be
(xy - yx) z = x(yz) - y (xz), 
\ee
or
\be
(xy) z - x (yz) = (yx) z - y (xz), \qquad  \forall \; x, y, z \in {\mathcal G}.
\ee
Thus, ${\mathcal G}$ is a {\it quasiassociative  algebra} , and $T^* {\mathcal G}: = {\mathcal G}
 \ {\mathop{\ltimes}\limits_\rho} \  {\mathcal G}^*$  (A2.1) is a 
proper phase space of~${\mathcal G}$, meaning that the symplectic form $\omega$ is a 2-cocycle 
on $T^* {\mathcal G}$.  This explains the appearance of quasiassociative algebras in~[7,~9].  
There are many examples of quasiassociative algebras given in~[7,~9] (in addition 
to obvious ones coming from associative algebras) such as Lie algebras of vector f\/ields on 
${\mathbf R}^n$~[7] and on $GL (n)$~[9].  We shall now use the former example 
${\mathcal D}_n = {\mathcal D} ({\mathbf R}^n) = 
\{X \in R^n, R = C^\infty ({\mathbf R}^n) \}$, with the quasiassociative multiplication
\be
(\pbf{X} \pbf{Y})^i = \sum_s X^s Y,^i_s 
\ee
where
\be
( \ \cdot \ ),_s = \partial_s ( \ \cdot \ ) = {\partial ( \ \cdot \ ) \over \partial x^s} . 
\ee
The quasiassociative property (A2.10) is satisf\/ied because 
\be
(\pbf{X} (\pbf{Y} \pbf{Z}) - (\pbf{X} \pbf{Y}) \pbf{Z})^i = \sum_{st} X^s Y^t Z,^i_{st}
\ee
is symmetric in $X$, $Y$.

\newpage

Now, formula (A2.11) means that
\be
\rho^d (\pbf{X}) = \hat X {\mathbf 1}, \qquad  \hat X = \sum_s X^s \partial_s.
\ee
Therefore, $\rho (X) = - \rho^d (X)^\dagger$ is: 
\be
 \rho(X) = \left(\sum_s \partial_s X^s\right) {\mathbf 1}.
\ee
Hence, the Lie bracket on the Lie algebra ${\mathcal D}^{(1)}_n = {\mathcal D}_n 
\ltimes {\mathcal D}_n^*$ is 
\be
\left[ {\pbf{X}  \choose \pbf{u}}, {\pbf{Y} \choose \pbf{v}} \right]^i_i = \sum_s 
\left(\matrix{X^s Y,^i_s  - Y^s X,^i_s \cr (X^s v_i - Y^s u_i),_s \cr} \right). 
\ee

By general theory, since the symplectic form $\omega$ (A2.2) is a  {\it nondegenerate} 
2-cocycle on ${\mathcal G}^{(1)} = {\mathcal G} \ltimes {\mathcal G}^*$, it can 
be represented by an ${\mathcal O}$-opertor 
${\mathcal O}:  {\mathcal G}^{(1)*} \rightarrow {\mathcal G}^{(1)}$.  
Writing elements of ${\mathcal G}^{(1)}$ as ${\alpha \choose a}$, 
$\alpha \in {\mathcal G}^*$, $a \in {\mathcal G}$, we see from formulae (A2.2) and (2.3) that 
\be
- \omega\left({X \choose u}, {Y \choose v} \right) =
 \langle  {\mathcal O}^{-1} {X \choose u}, {Y \choose v} \rangle  
= - \langle u, Y\rangle  + \langle v, X\rangle , 
\ee
whence
\be
{\mathcal O}^{-1} {X \choose u} = {-u \choose X}, 
\ee
so that
\be
{\mathcal O}^{-1} = \left(\matrix{{\mathbf{0}} & - {\mathbf{1}} \cr {\mathbf{1}} & 
{\mathbf{0}} \cr} \right), 
\ee
and hence
\be
{\mathcal O} = \left(\matrix{{\mathbf{0}} & {\mathbf{1}} \cr -{\mathbf{1}} & {\mathbf{0}} \cr} \right).
 \ee
The resulting Lie algebra bracket (2.11) on ${\mathcal G}^{(1)^*}$ takes the form
\be
\ba{l}
\ds \left[ {\alpha \choose a} , {\beta \choose b} \right] = {a \choose - \alpha}^\cdot 
{\beta \choose b} - {b \choose - \beta}^\cdot {\alpha \choose a} 
\vspace{3mm}\\
\ds \qquad {\mathop{=}\limits^{\mbox{\scriptsize [by (A2.23)]}}} \
\left(\matrix{a \beta - \beta a - \alpha b \cr a b \cr} \right) - 
\left(\matrix{b \alpha - \alpha b - \beta a \cr b a \cr} \right) = \left(\matrix{a \beta - b \alpha 
\cr [a, b] \cr} \right).
\ea
\ee

\noindent
{\bf Lemma A2.22.} {\it  In ${\mathcal G}^{(1)*}$, the coadjoint action formulae is
\setcounter{equation}{22}
\be
{X \choose u}^\cdot {\alpha \choose a} = {X \alpha - \alpha X + ua \choose Xa}, \qquad
X, a \in {\mathcal G}, \quad  u, \alpha \in {\mathcal G}^*.
\ee
(The new notation is explained in formulae (A2.25) below.) }

\medskip

\noindent
{\bf Proof.} We have,
\be
\ba{l}
\ds \langle  {X \choose u}^\cdot {\alpha \choose a}_, {Y \choose v} \rangle 
 \sim  - \langle  {\alpha \choose a}, 
\left[ {X \choose u}_, {Y \choose v} \right] \rangle  
= - \langle  {\alpha \choose a}, {[X, Y] \choose 
\rho (X) (v) - \rho (Y) (u)} \rangle  
\vspace{3mm}\\
\qquad = - \langle  \alpha, [X, Y] \rangle  - \langle  \rho (X) (v) - \rho (Y) (u), a \rangle   
\vspace{3mm}\\
\qquad \sim  \langle  X^. \alpha, Y\rangle   - \langle v, \rho (X)^\dagger (a) \rangle  
+ \langle  u, \rho (Y)^\dagger (a) \rangle  
 = \langle  X^. \alpha, Y\rangle  + \langle  v, X a\rangle  - \langle u, Ya \rangle . 
\ea\hspace{-22.45pt}
\ee
Let us def\/ine the left and right multiplication of ${\mathcal G}$ on ${\mathcal G}^*$ 
by the relations
\renewcommand{\theequation}{{\rm A2}.\arabic{equation}{\rm a}}
\setcounter{equation}{24}
\be
\langle  X u, Y \rangle   \sim  - \langle  u, XY \rangle , \qquad \forall \;
u \in {\mathcal G}^*, \  X, Y \in {\mathcal G}, 
\ee
\renewcommand{\theequation}{{\rm A2}.\arabic{equation}{\rm b}}
\setcounter{equation}{25}
\be
\langle  u X, Y \rangle  \sim  \langle  u, Y X \rangle , \qquad
 \forall\;  u \in {\mathcal G}^*, \ X, Y \in {\mathcal G}. 
\ee
Then
\renewcommand{\theequation}{{\rm A2}.\arabic{equation}}
\setcounter{equation}{25}
\be
\langle  X^. u, Y\rangle   \sim 
 \langle  u, [Y, X] \rangle  = \langle  u, YX - XY \rangle   \sim  \langle  - u X + Xu, Y \rangle , 
\ee
so that
\be
X^. u = X u - uX. 
\ee
Substituting formulae (A2.25b, (A2.27) into formula (A2.24), we get
\be
\langle  X \alpha - \alpha X + ua, Y\rangle  + \langle  v, X a \rangle , 
\ee
and formula (A2.23) results. \hfill \rule{3mm}{3mm}

Now, formula (A2.25a) yields
\be
\ba{l}
\langle  Xu, Y \rangle  \sim  - \langle  u, XY\rangle  
= - \langle  u, \rho^d (X) (Y) \rangle 
\vspace{2mm}\\
\qquad  = - \langle  u, - \rho (X)^\dagger (Y) \rangle  
\sim  \langle  \rho (X) (u), Y \rangle  , 
\ea
\ee
so that
\be
Xu = \rho (X) (u). 
\ee
Formula (A2.21) can be now rewritten as
\be
\left[{\alpha \choose a}, {\beta \choose b}\right] =  {\rho(a) (\beta) - 
\rho (b) (\alpha) \choose [a, b]}.
\ee
We see that the Lie algebra structure on ${\mathcal G}^{(1)*}$
is the same as on ${\mathcal G}^{(1)}$.  This 
explains, -- and  provides an $n$-dimensional generalization of, -- formula (3.17) for the case 
$\epsilon = \mu = 0$.  (Formula (A2.16) is an $n$-dimensional analog of formula~(3.7).)

We conclude by def\/ining the notion of a {\it symplectic  double} of the Lie algebra 
${\mathcal G}^{(1)} = {\mathcal G} \ltimes {\mathcal G}^*$ for a quasiassociative ${\mathcal G}$. 
It is dif\/ferent from the Drinfel'd's 
classical double, ${\mathcal G} + {\mathcal G}^*$, discussed in the previous Section.

So, let ${\mathcal A}$ be a quasiassociative ring and ${\mathcal G} = \mbox{Lie}\, (A)$. 
 Def\/ine $T^* {\mathcal A}$ by the formula~[7]
\be
{a \choose a^*} {b \choose b^*} = {ab \choose ab^*}, \qquad
a, b \in {\mathcal A}, \quad  a^*, b^* \in {\mathcal A}^*,
\ee
where the product $ab^*$ is def\/ined by formula (A2.25a).  Notice, that $a^* b = 0$ in 
formula (A2.32), in contrast to the previous def\/inition (A2.25b).

\medskip

\noindent
{\bf Proposition A2.33.} {\it $T^* {\mathcal A}$  (A2.32) is again a quasiassociative algebra. }

\medskip

\noindent
{\bf Proof.}  We have:
\setcounter{equation}{33}
\be
\ba{l}
\ds \left({a \choose a^*} {b \choose b^*} \right) {c \choose c^*} - {a \choose a^*} \left( 
{b \choose b^*} {c \choose c^*} \right) 
\vspace{3mm}\\
\ds \qquad = {ab \choose ab^*} {c \choose c^*} - {a \choose a^*} {b c \choose b c^*} = 
{(ab) c - a (bc) \choose (ab) c^* - a (bc^*)} .
\ea
\ee
Thus, we need to show that
\be
(ab)c^* - a (bc^*) = (ba) c^* - b (ac^*). 
\ee
We have, $\forall \; d \in {\mathcal A}$:
\be
\langle (ab)c^* - a(bc^*), d \rangle   \sim 
 \langle c^*, - (ab)d \rangle  + \langle bc^*, ad \rangle  
\sim \langle c^*, - (ab) d - b (ad) \rangle  . 
\ee
But
\be
(ab) d + b (ad) = (ba) d + a (bd)
\ee
by formula (A2.9). \hfill \rule{3mm}{3mm}

Since $\mbox{Lie}\,(T^* {\mathcal A}) = \mbox{Lie}\, ({\mathcal A}) \ltimes [\mbox{Lie}\,
 ({\mathcal A})]^*$:
\be
\left[{a \choose a^*}, {b \choose b^*} \right] = {a \choose a^*} {b \choose b^*} - 
{b \choose b^*} {a \choose a^*} = {ab - ba \choose ab^* - ba^*}, 
\ee
we can construct a symplectic 2-cocycle, -- and thus an ${\mathcal O}$-opertor, -- starting with 
${\mathcal G}^{(1)} = {\mathcal G} \ltimes {\mathcal G}^*$ rather than ${\mathcal G}$.  
This process can be continued indef\/initely.  
It is natural to call the Lie algebra ${\mathcal G}^{(1)} \ltimes {\mathcal G}^{(1)^*}$ 
the {\it symplectic double} of the Lie algebra ${\mathcal G}^{(1)} = {\mathcal G} \ltimes 
{\mathcal G}^*$; thus, ${\mathcal G}^{(1)}$ is the 
symplectic double of ${\mathcal G}$.  According to results of \S~6, the space
$\mbox{Fun}\,({\mathcal G}^{(1)*})$ 
carries three compatible Hamiltonian structures:  symplectic,  linear, and quadratic.

\label{kupershmidt_7-lp}

\end{document}